\documentclass{amsproc}
\usepackage{amssymb}
\usepackage{latexsym, amsmath, amscd,amsthm}
\usepackage{graphicx}
\usepackage[percent]{overpic}

\newtheorem{theorem}{Theorem}
\newtheorem*{maintheorem}{Main Theorem}
\newtheorem{lemma}[theorem]{Lemma}
\newtheorem{proposition}[theorem]{Proposition}
\newtheorem{corollary}[theorem]{Corollary}

\theoremstyle{definition}
\newtheorem*{definition}{Definition}

\theoremstyle{remark}

\def\defn#1{Definition~\ref{def:#1}}
\def\thm#1{Theorem~\ref{thm:#1}}
\def\lem#1{Lemma~\ref{lem:#1}}
\def\figr#1{Figure~\ref{fig:#1}}
\def\prop#1{Proposition~\ref{prop:#1}}
\def\cor#1{Corollary~\ref{cor:#1}}
\def\sect#1{Section~\ref{sect:#1}}

\newcommand{\R}{\mathbb{R}}
\newcommand{\arc}[1]{\gamma_{#1}}

\newcommand{\bdry}{\partial}

\newcommand{\setm}{\smallsetminus}

\newcommand{\m}{\mathcal}

\let\mgp=\marginpar \marginparwidth18mm \marginparsep1mm
\def\marginpar#1{\mgp{\raggedright\tiny #1}}

\makeatletter
\newbox\overbox
\def\fakeover#1{\setbox\overbox\hbox{$#1$}\hbox
                         {$\overline{#1\hskip-\wd\overbox}$\hskip\wd\overbox}}
\def\overnoarrow#1{\mskip1.5mu\overline{\mskip-1.5mu#1}}
\def\overrightarrow#1{\mskip2mu\vbox{\m@th\ialign{##\crcr
   \rightarrowfill\crcr
   \noalign{\kern-.4pt               
        \kern-\fontdimen22\textfont2 
        \nointerlineskip}
   ${\mskip0mu\hfil\fakeover{#1}\hfil\mskip6mu}$\crcr}}\mskip-2mu}
\def\overleftrightarrow#1{\mskip-3mu\vbox{\m@th\ialign{##\crcr
   \leftarrowfill\hskip-.6em\rightarrowfill\crcr
   \noalign{\kern-.4pt               
        \kern-\fontdimen22\textfont2 
        \nointerlineskip}
   ${\mskip3mu\hfil\fakeover{#1}\hfil\mskip6mu}$\crcr}}\mskip-2mu}

\def\segment#1{\smash{\overnoarrow{#1}}}
\makeatother

\newcommand{\altval}{15.66}

\begin{document}
\title[Alternating quadrisecants - preprint]
     {Alternating Quadrisecants of Knots }

\author{Elizabeth Denne}
\address{Mathematics Department \\
         Harvard University \\
         Cambridge, MA, 02138 }
\email{denne@math.harvard.edu}

\subjclass{Primary 57M25.}

\keywords{Knots, links, quadrisecants, total curvature}

\begin{abstract} It is known \cite{Pann,Kup} that for every knotted
  curve in space, there is a line intersecting it in four places, a
  quadrisecant. Comparing the order of the four points along the line
  and the knot we can distinguish three types of quadrisecants; the
  alternating ones have the most relevance for the geometry of a knot.
  I show that every (nontrivial, tame) knot in $\R^3$ has an
alternating quadrisecant. This result has applications to the total
  curvature, second hull and ropelength of knots.   
\end{abstract}
\date{\today}
\maketitle

\section{Introduction}\label{sect:intro} 

Throughout this paper a \emph{knotted curve} will denote an oriented
nontrivial tame knot in $\mathbb{R}^3$. (By knot, we mean a
homeomorphic image of $S^1$ in $\R^3$, modulo reparametrizations. By
tame, we mean the knot is ambient isotopic to a polygonal knot.)  A
secant line is a straight line which intersects the knot in at least
two distinct places. Trisecant, quadrisecant and quintisecant lines
are straight lines which intersect a knot in at least three, four and
five distinct places, respectively. It is clear that any closed curve
has a 2-parameter family secants. A simple dimension count shows there
is a 1-parameter family of trisecants and that
quadrisecants are discrete (0-parameter family). Quintisecants are a
$-1$-parameter family; meaning they exist only for a codimension~1 set
of knots (\prop{pfgeneric}). 

Thus we expect knotted curves to have quadrisecants. Indeed,
in 1933, E.~Pannwitz~\cite{Pann} first showed that every nontrivial
generic polygonal 
knot in~$\mathbb{R}^3$ has at least $2u^2$ quadrisecants, where
$u$ is the unknotting number of a knot.  In~1980, H.R.~Morton and
D.M.Q.~Mond~\cite{MM} independently proved every nontrivial generic knot
has a quadrisecant and conjectured that a generic 
knot with crossing number~$n$ has at least~${n}\choose{2}$
quadrisecants. It was not until 1994 that G.~Kuperberg~\cite{Kup} managed
to extend the result and show that all (nontrivial tame) knots in
$\mathbb{R}^3$ have a quadrisecant. 

Quadrisecants come in three basic types. These are distinguished
by comparing the orders of the four points along the knot and along the
quadrisecant line. There are three relative
orderings of an oriented quadrisecant and an
unoriented knot. If $abcd$ is the ordering along the quadrisecant,
then the possible orderings along the knot are $abcd$, $abdc$ and
$acbd$. I will show that
every nontrivial tame knot has at 
least one quadrisecant with knot ordering $acbd$. This type of
quadrisecant is called an alternating quadrisecant. This result
refines the previous work 
about quadrisecants, giving greater geometric insight into knots and
also providing several interesting applications. This paper contains the main results of my PhD thesis \cite{Den} and I am grateful to my advisor J.M~Sullivan for suggesting this problem to me.  Recently, R.~Budney~\emph{et
al}.~\cite{BCSS} have shown that the finite type $2$ Vassiliev
invariant can be computed by counting alternating
quadrisecants with appropriate multiplicity. While this result implies the
existence of alternating quadrisecants for many knots, our Main
Theorem shows existence for \emph{all} (nontrivial tame) knots. On
the other hand, our results provide no way to count alternating
quadrisecants.  

The existence of alternating quadrisecants provides new proofs to two
previously known results about the geometry of knotted curves.
For smooth closed curves, the total curvature can be thought of as the
total angle through which the unit tangent vector turns (or the length
of the tangent indicatrix). Equivalently, the total curvature of an arbitrary curve is the supremal total curvature of inscribed polygons. Around~1949 I.~F\'ary~\cite{Fary}
and J.W.~Milnor~\cite{Mil} proved 
independently that the total curvature of a nontrivial tame knot in
$\mathbb{R}^3$ is at least~$4\pi$. Note that the total
curvature of an inscribed polygon is 
less than or equal to the total curvature of the curve it is inscribed
in. The total curvature of an alternating
quadrisecant is~$4\pi$. Thus if a knotted curve has an alternating
quadrisecant, the total curvature will be greater than or equal
to~$4\pi$, giving the F\'ary-Milnor theorem another proof.
In~1998, the F\'ary-Milnor theorem was extended to knotted curves in Hadamard\footnote
{A Hadamard manifold is a complete simply-connected
Riemannian manifold with non-positive sectional curvature.}
manifolds by C.~Schmitz~\cite{Schm} and S.B.~Alexander and
R.L.~Bishop~\cite{AB}. 

The second new proof is for the theorem about the second hull of a
knotted curve in $\mathbb{R}^3$. Intuitively, the second hull is the
part of space that the knotted curve winds around twice. In 2000,
J.~Cantarella \emph{et al}~\cite{ckks} proved that the second hull of
a knotted curve in $\mathbb{R}^3$ is nonempty. This paper conjectured
the existence of alternating quadrisecants for knotted curves in
$\mathbb{R}^3$ as another way of proving that the second hull is
nonempty. This is because the mid-segment $bc$ of alternating
quadrisecant $abcd$ is in the second hull.

Finally, quadrisecants may be applied to the ropelength
problem. The ropelength problem asks to minimize the length of a
knotted curve subject to maintaining an embedded tube of fixed diameter
around the tube; this is a mathematical model of tying the knot tight
in a rope of fixed thickness. The ropelength of a knot is the quotient 
of its length and its thickness. The thickness is the diameter of the
largest embedded normal tube around the knot. The exact value of
the minimum ropelength needed to tie any nontrivial knot is not
currently known. In joint work with J.M.~Sullivan and Y.~Diao~\cite{dds},
essential alternating quadrisecants are used to improve
the known lower bounds of ropelength from~$12$ to~\altval. Several
independent simulations (see for instance \cite{Pie,Sul})
have found a trefoil knot with ropelength less than~$16.374$. This is
presumably close to the minimizer, so the new bounds are quite sharp.

The proof of the main theorem will extend ideas taken from all of the previous
papers, but in particular from~\cite{Kup, Pann, Schm}. At its core, the
proof assumes that alternating quadrisecants do not exist, then uses
this to create a contradiction to knottedness. A quadrisecant includes
a number of trisecants, thus a large part of the 
proof will be dedicated to a detailed understanding the structure of
the set of trisecants of a knot, both in $K^3$ and when
projected to the set of secants~$S=K^2\setm{\Delta}$. 

\sect{secants} introduces terminology and explores the relationship
between quadrisecants and trisecants. \lem{tri1} shows any knotted
curve has at least a 1-parameter family of trisecants
and \lem{alt} shows that alternating
quadrisecants occur when trisecants of same and different
orders share common points. A detailed understanding of the structure
of the set of trisecants is required. In 
particular, this set has very nice properties when the knot is a
generic polygonal knot. In \sect{genpoly}, the definition of a generic
polygonal knot is given and the conditions for genericity are
shown to be generic. In \sect{tri} the structure of the set of trisecants of
generic polygonal knots is examined in detail, both in
$K^3$ and $K^2\setm {\Delta}$. In $K^2\setm {\Delta}$, it is shown to be a piecewise immersed 1-manifold that
intersects itself in double points. In the end we wish to
prove the existence of alternating quadrisecants for any nontrivial
tame knot. Any tame knot is the limit of a sequence of generic
polygonal knots. Thus we first prove the
existence of an essential alternating quadrisecant for nontrivial generic
polygonal knots. The notion of essential is required so that
quadrisecants do not degenerate to trisecants in the
limit. \sect{equad} defines the notion of essential for secants, 
trisecants and quadrisecants and describes the structure of the set of
essential trisecants.  \sect{quadgenpoly}
gives the proof that any nontrivial generic polygonal knotted curve has an
essential alternating quadrisecant. In \sect{mainresult}, the Main Theorem
(\thm{mainthma}) is proved: every nontrivial tame knot in~$\R^3$ has
an alternating quadrisecant. This is strengthened in \cor{altquadess}: every
 nontrivial knot of finite total curvature in~$\R^3$ has an
\emph{essential} alternating quadrisecant. Finally,
\sect{corollaries} provides 
applications of essential alternating quadrisecants to the total
curvature, second hull and ropelength of a knotted curves. 


\section{Secants, trisecants and quadrisecants}\label{sect:secants}

\begin{figure}
\begin{overpic}[scale=.4]{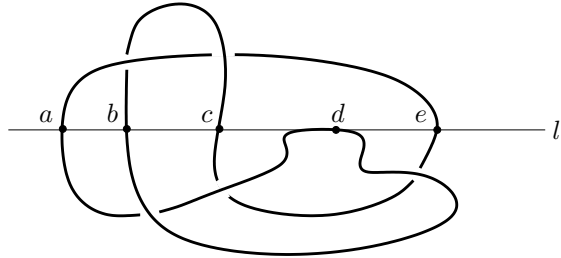}
\put(6,25){$a$}
\put(18.5,25){$b$}
\put(36,25){$c$}
\put(60,25){$d$}
\put(75.5,25){$e$}
\put(101,22){$l$}
\end{overpic}
\caption{A quintisecant line $l$ and quintisecant $abcde$.}
\label{fig:secantline}
\end{figure}

Recall that a \emph{knotted curve} is an
oriented nontrivial tame knot in~$\mathbb{R}^3$.
 
\begin{definition} An \emph{$n$-secant line} for a knotted curve $K$ is an
oriented line whose intersection with $K$ has at
least $n$ components.
\end{definition}

\begin{definition} An  \emph{$n$-secant} is an ordered $n$-tuple of
  points in $K$ (no two of which lie in a common straight subarc of
  $K$) which lie in order on an $n$-secant line.
A $2$-secant will be called a secant, a $3$-secant will be called a
trisecant, a $4$-secant will be called a quadrisecant and a $5$-secant
will be called a quintisecant.
\end{definition}

The line through an $n$-secant may intersect the knot in more than $n$
places, but this does not affect the definition of an $n$-secant.
\figr{secantline} shows a knot with quintisecant line~$l$
and quintisecant $abcde$. For example, the quintisecant line $l$ also
includes quadrisecant~$abcd$, trisecant~$abe$, and secant~$bd$. 

In $K^n$, let $\tilde{\Delta}$ denote the set of $n$-tuples in which some
pair of points lie in a common straight subarc of~$K$. Also let
$\Delta$ denote the big diagonal, the set of $n$-tuples in which some
pair of points are equal. Then $\Delta\subset \tilde{\Delta}$. (The
inclusion can be strict as
for polygonal knots.) For a closed curve~$K$, any two distinct
points determine a straight line. These points are a secant if and
only if they do not lie on a
common straight subarc of~$K$. Thus the set of
secants $S=K^2\setm\tilde{\Delta}$ and is 
topologically an annulus. More generally, any $n$~distinct
points are not necessarily collinear, hence the set of $n$-secants
is contained in~$K^n\setm\tilde{\Delta}$. 

\begin{figure}
\begin{overpic}[scale=.4]{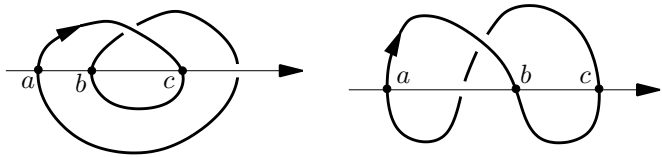}
\put(2.5,10.2){$a$}
\put(10.7,9.7){$b$}
\put(24,10.3){$c$}
\put(59.3,11){$a$}
\put(78,11){$b$}
\put(87,11){$c$}
\end{overpic}
\caption[Trisecants of same and different ordering.]{The points of each
  trisecant have order $abc$ along the line. The left has
  \emph{different} order $acb$ along the knot and the right has
  \emph{same} order $abc$ along the knot.}
\label{fig:trisecants} 
\end{figure}

Consider the set of trisecants. Take a trisecant $abc$ with
intersection points in that order along the trisecant line. There are 
$\vert S_3\vert =6$ possibilities for their
ordering along~$K$. Along~$K$, the order is just a cyclic order. Thus
there are~$\vert S_3/C_3\vert=2$ cyclic orderings of the oriented knot
and oriented trisecant line. Picking the lexicographically least element
in each coset, we see the order along the knot is~$abc$ or~$acb$. These orderings are respectively called {\it same} and {\it
different}. \figr{trisecants} illustrates the two types of
trisecant. 

We may make a more general observation about the ordering of points
along a knot. Let $abc\in  K^3$ and assume that
$a$, $b$ and~$c$ are distinct. As~$K$ is oriented, moving from~$a$ to~$b$ to~$c$ along the $K$ will either match the orientation of~$K$ or
not. These are the two cyclic orderings discussed above. Let~$\mathcal{S}$ be the set of triples of~$K^3$ where the
triples have the {\it same} ordering along the knot and let~$\mathcal{D}$ be
the set of triples of~$K^3$ where the ordering of the
triples {\it differs} from their ordering along the knot. Observe that
$\mathcal{S}$ and $\mathcal{D}$ {\it are} the connected components of~$K^3\setm\tilde{\Delta}$.  Then~$\mathcal{S}$ and~$\mathcal{D}$ are
disconnected as perturbing~$a$,~$b$ and~$c$ a
little will not change their ordering. The only way to change the
order along the knot from $abc$ to $acb$ is for $a$, $b$ and $c$ to no
longer be distinct. That is, either $a=b$, $b=c$ or $a=c$. This is
precisely the big diagonal~$\Delta\subset\tilde{\Delta}$.

\begin{definition}
Let $\mathcal{T}\subset K^3\setm\tilde{\Delta}$ denote the set of all
  trisecants of a knot~$K$. Let the closure of~$\mathcal{T}$ in~$K^3$ be denoted by~$\overline{\mathcal{T}}$ and let the boundary of $\mathcal{T}$ be defined as $\bdry \mathcal{T} := \overline{\mathcal{T}}\setm\mathcal{T}$. Let
 $\mathcal{T}^s=\mathcal{T}\cap\mathcal{S}$ denote the set of
  trisecants where the ordering along the trisecant is the \emph{same}
  as the ordering along the knot.  Let
  $\mathcal{T}^d=\mathcal{T}\cap\mathcal{D}$ denote the set of all
  trisecants where the ordering along the trisecant is \emph{different}
  to the ordering along the knot. Clearly $\mathcal{T}^s\cap
  \mathcal{T}^d =\emptyset$. Note that switching the orientation of
  either the knotted curve or the trisecant line interchanges~$\mathcal{T}^s$ and~$\mathcal{T}^d$.
\end{definition}

\begin{figure}
\begin{overpic}[scale=.5]{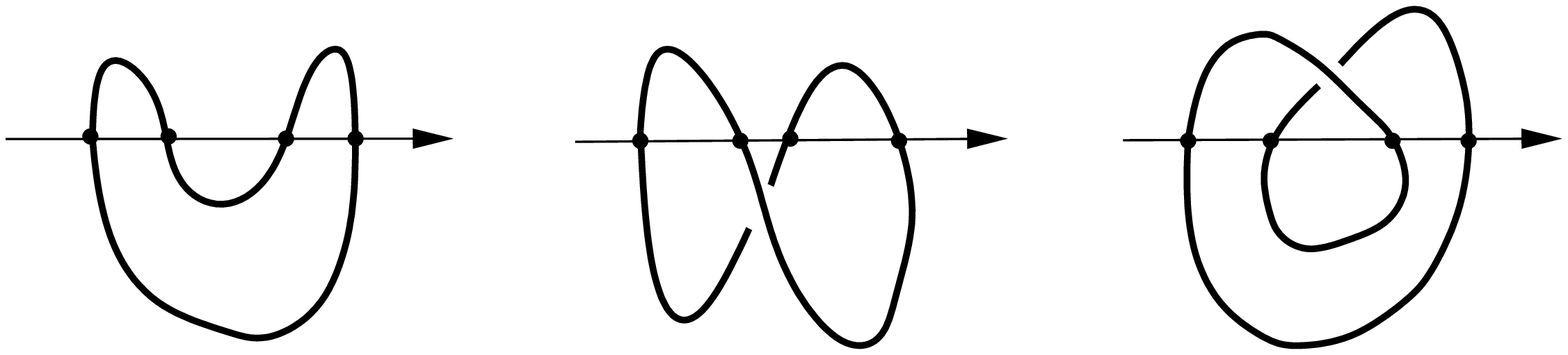}
\put(3.3,11){$a$}
\put(8.5,10.7){$b$}
\put(18.5,11){$c$}
\put(23.5,10.7){$d$}
\put(38.5,11){$a$}
\put(45.4,10.7){$b$}
\put(50.5,11){$c$}
\put(55.3,10.7){$d$}
\put(73.5,11){$a$}
\put(78.8,10.7){$b$}
\put(87.5,11){$c$}
\put(91.5,10.7){$d$}
\end{overpic}
\caption[Three types of quadrisecants: simple, flipped and alternating.]{From left to right, quadrisecants~$abcd$ are simple, flipped and alternating,
because the dihedral order of the points along~$K$ is~$abcd$,~$abdc$
and~$acbd$, respectively.}
\label{fig:quadord}
\end{figure}

Just as with trisecants, we may compare the ordering of the
intersection points of a quadrisecant line with their
ordering along the knot. Given intersection points~$abcd$
in that order along the quadrisecant line, there are~$|S_4|$
possibilities for their ordering along~$K$. Again, the order along~$K$
is only a cyclic order and ignoring the orientation of~$K$ is just a
dihedral order. Thus there are~$\left|S_4/D_4\right|=3$ dihedral
orderings of an oriented quadrisecant and unoriented
knot. Let~$abcd$ 
be the order of intersection points along the quadrisecant line, then
the three equivalence classes or types of quadrisecants are
represented by $abcd$, $abdc$ and $acbd$, where we have chosen the
lexicographically least order as the name for each.  \figr{quadord}
illustrates these orderings.  

\begin{definition} Quadrisecants of type~$acbd$ are called \emph{alternating quadrisecants}. Quadrisecants of types~$abcd$ and~$abdc$ are called \emph{simple} and \emph{flipped} respectively.
\end{definition}

 Alternating quadrisecants
(see the right-most quadrisecant in \figr{quadord}) are so named as the 
quadrisecant's points alternate from one end of the quadrisecant line to
the other as the ordering along the knot is followed. Note that unlike
trisecants, switching the orientation of the quadrisecant line does not
change the type of quadrisecant.

We have already noted that each closed curve has a 2-parameter family of secants~$S=K^2\setm\tilde{\Delta}$. It turns out that knotted curves also have
(at least) a 1-parameter family of trisecants as seen in the following
lemma originally due to Pannwitz~\cite{Pann}. 

\begin{lemma}\label{lem:tri1} Each point of a knotted curve~$K$ is the 
first point of at least one trisecant. 
\end{lemma}
\begin{proof} Suppose there is a point~$a$ of the knot which is not
the start point of any trisecant. That is, no points~$b, c\in K$ are
collinear with~$a$, with~$a$ as the start point. The union of all
chords~${ab}$ for~$b\in K$ is a disk with boundary~$K$. If two chords~$ab$ and~$ac$ intersect at a place other than~$a$, then they overlap
and one is a subinterval of another. They form a trisecant ($abc$ or~$acb$), contrary to the assumption. Thus the disk is embedded,
contradicting the assumption of knottedness.
\end{proof}

If we view the knotted curve~$K$ from a point~$a\in K$ as in
\lem{tri1}, then the trisecants occur where we see a 
crossing of~$K$. Thus we expect that there will be a finite set of
trisecants with first point~$a$. In fact, Pannwitz~\cite{Pann} proved for
generic polygonal knots that there are at least $2u^2$ trisecants from
any point~$a\in K$ (where~$u$ is the unknotting number of the
knot). Schmitz~\cite{Schm} has shown
\lem{tri1} also holds for knotted curves in Hadamard manifolds. In
\sect{tri} we will show that the set of trisecants is generically an embedded
1-manifold in~$K^3$.  

\begin{figure}
\begin{overpic}[scale=.4]{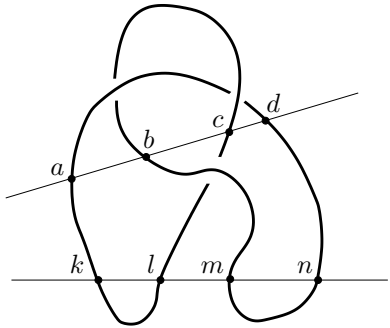}
\put(12,39.5){$a$}
\put(36,46.5){$b$}
\put(54,52){$c$}
\put(68,56){$d$}
\put(17,14){$k$}
\put(37,14){$l$}
\put(51,14){$m$}
\put(76,14){$n$}
\end{overpic}
\caption[Non-alternating and alternating quadrisecants.]{A trefoil
  knot with  alternating quadrisecant~$abcd$ and simple
  quadrisecant~$klmn$. Other trefoil knots might not have the simple
  quadrisecant 
  shown. However, the Main Theorem shows that every knotted curve has
  an alternating quadrisecant.}\label{fig:twoquad}
\end{figure}

What about higher order secants? In \sect{genpoly} we will see that
knotted curves generically have no $n$-secants where $n\geq
5$. However, knotted curves do have quadrisecants.
\figr{twoquad} shows a trefoil with two quadrisecants.
Quadrisecant~$abcd$ is an alternating quadrisecant and
quadrisecant~$klmn$ is a simple
quadrisecant. The trefoil knot, for example, does not necessarily
have the simple quadrisecant shown. The rest of this work will be
dedicated to proving the following theorem: 

\begin{maintheorem} Every knotted curve in~$\mathbb{R}^3$
has at least one alternating quadrisecant. 
\end{maintheorem}

Observe that quadrisecants are formed when several trisecants share common
points. Quadrisecant~$abcd$ includes four trisecants (1)$abc$,
(2)$abd$, (3)$acd$, (4)$bcd$. Pannwitz~\cite{Pann} showed
quadrisecants exist by looking for pairs of trisecants like~(1)$abc$
and~(3)$acd$. Here, the first and third points of trisecant~$abc$ are
the same as the first and second points of trisecant~$acd$. Kuperberg~\cite{Kup} showed that quadrisecants exist by looking for pairs
of trisecants like~(2)$abd$ and~(3)$acd$. Here, the first and third
points of the trisecants are the same. Schmitz~\cite{Schm} nearly
showed that alternating quadrisecants exist by looking at families~(1)$abc$ and~(2)$abd$, where the first and second points of the trisecants are the
same, but in his proof, some quadrisecants may degenerate to trisecants. 
I will use Schmitz' approach, as it allows us to use the
orderings of coincident trisecants to determine the ordering of the
quadrisecant.

Let quadrisecant~$abcd$ be an alternating quadrisecant as in
\figr{quadord} (right). As we have seen it includes several trisecants: ~$abc$ and~$bcd$ are trisecants of different ordering, and~$abd$ and~$acd$ are trisecants of same ordering. Given an alternating
quadrisecant, there will be trisecants of same and different ordering
which have common points. For example,~$abd\in\mathcal{T}^s$ and~$abc
\in\mathcal{T}^d$ have the first two points in common. In fact the
converse is true: if trisecants of same and different orderings have
the first two points in common, then there is an alternating quadrisecant.

Define the projection $\pi_{12}: K^3 \rightarrow K^2$ by
$\pi_{12}(xyz)=xy$ and let  $T=\pi_{12}(\mathcal{T})\subset S$ be the
projected image of trisecants~$\mathcal{T}$. A secant~$ab\in S$ is in~$T$ if and only if it is the first two points of some trisecant~$abc$.
We further set~$T^s:=\pi_{12}(\mathcal{T}^s)$ and~$T^d:=\pi_{12}(\mathcal{T}^d)$.

\begin{lemma}\label{lem:alt} Let $ab\in T^s\cap T^d$ in~$S$. This
  means that there exists~$c$,~$d$ such that~$abc\in \mathcal{T}^d$ and~$abd\in\mathcal{T}^s$. Then either~$abcd$ or~$abdc$ is an alternating
quadrisecant. 
\end{lemma}
\proof Trisecants $abc$ and $abd$ lie on a common line, which must be
a quadrisecant line. (Points~$c$ and~$d$ do not lie on the same
straight subarc of~$K$. If they did, then trisecants~$abc$ and~$abd$
would have the same ordering.) The order of the intersection points
along the quadrisecant line is either~$abcd$ or~$abdc$.  First, assume the order is~$abcd$. Using the definition of~$T^s$ and~$T^d$, the order of the
intersection points along the knot must be $acbd$. But this means that~$abcd$ is an alternating quadrisecant. A similar argument shows~$abdc$ is an alternating quadrisecant.  
\endproof

Thus to prove that any knotted curve has at least one alternating
quadrisecant, it is sufficient to prove~$T^s\cap T^d \neq\emptyset$ in~$S$. To do this, we need to restrict our attention to generic
polygonal knots (to be 
defined in the next section). We will then prove the existence of
alternating quadrisecants for generic polygonal knots and will use limit
arguments to extend the result to all knotted curves.

\section{Generic polygonal knots}\label{sect:genpoly}

In order to clearly understand the structure of the set of trisecants
we only consider generic polygonal
knots. This section defines generic polygonal knots and looks at some
implications of the definition for $n$-secants. The next section uses the
conditions which define generic polygonal knots to help give a
detailed description of the structure of the set of trisecants. 
 
\begin{definition} Let $K$ be a closed polygonal curve in
  $\mathbb{R}^3$. Let the vertices of~$K$ be
  $v_1,\dots,v_n\in\mathbb{R}^3$. Then~$K$ is \emph{non-degenerate} if no
  four~$v_i$ are coplanar and no three~$v_i$ are collinear. 
\end{definition}

An $n$-gon in $\R^3$ is determined by the position of its~$n$
vertices. Thus we identify the space of all $n$-gons with~$\R^{3n}$
with the usual topology. 

\begin{proposition}\label{prop:dense} The set of non-degenerate
  $n$-gons is an open dense set in~$\mathbb{R}^{3n}$.
\end{proposition}
 
\proof 
When $n=3$, the third vertex cannot lie on the straight line spanned
by the other two vertices. Thus the set of
degenerate configurations has codimension~1 in the set of all
configurations.

For $n\geq 4$, the set of degenerate $n$-gons is the union of~$n\choose 4$ cubic hypersurfaces, each given as the locus where some
particular quadruple of vertices is coplanar. This has codimension~1
in the set of all configurations. Note that the situation when three
vertices are 
collinear is a subset of the codimension~1 cubic hypersurfaces, so a special case is not needed to cover this situation. (They will have
even higher codimension.)  

Thus the set of all degenerate configurations is a codimension~1 (or higher)
algebraic surface in the set of all configurations. This is closed and nowhere
dense. Hence the non-degenerate $n$-gons are open and dense in~$\R^{3n}$.
\endproof

\begin{definition}
An edge $e_k$ denotes a \emph{closed} edge, including its
endpoints. One vertex is \emph{consecutive} to another vertex if they share a
common edge. An edge is \emph{adjacent} to another edge if they share a
common vertex. Thus given a vertex~$v$, there are two edges incident
to it (and adjacent to each other). The \emph{osculating plane} at a vertex~$v$ is the plane spanned by the two edges incident to~$v$. 
\end{definition}

\begin{proposition} Any non-degenerate polygon $K$ has the following
  properties: 
\item{(1)} It is embedded.
\item{(2)} The line connecting any two consecutive vertices only
  intersects~$K$ in their common edge. The line connecting any other
  two vertices does not intersect~$K$ again. 
\item{(3)} Two coplanar edges always share a common vertex.
\item{(4)} A line intersecting the interiors of two adjacent edges
  does not hit any vertex of~$K$.
\item{(5)} The three points of a trisecant do not lie on three adjacent edges.
\end{proposition}

\proof 
\item{(1)} If $K$ is not embedded, then~$K$ intersects itself. If
  two non-adjacent edges edges~$e_i$,~$e_j$ intersect, 
then the four vertices of~$e_i$ and~$e_j$ are coplanar
contradicting non-degeneracy. Note that two adjacent edges can not
  intersect (except at their common vertex) as they are not collinear.
\item{(2)} Any two consecutive vertices share a common edge~$e_k$
  which determines a line~$E_k$. Suppose~$E_k$ intersects~$K$ in
  another edge~$e_l$. Then the vertices of~$e_k$ and~$e_l$ are
  coplanar, contradicting non-degeneracy. Now suppose the line
  connecting any two (non consecutive) vertices~$v_i$ and~$v_j$
  intersects~$K$ in edge~$e_l$. Then~$v_i$,~$v_j$ and the vertices of~$e_l$ are coplanar, again contradicting non-degeneracy. 
\item{(3)} Two coplanar edges always share a
common vertex, else the plane contains four vertices of a
non-degenerate polygon. In other words, the lines determined by any two non-adjacent edges of a non-degenerate polygon are skew.
\item{(4)} A line intersecting the interiors of two adjacent edges
  does not hit any vertex of~$K$. If it did then the plane containing
  the edges also contains the vertex. 
\item{(5)} The three points of a trisecant do not lie on three adjacent edges. Suppose trisecant~$t$ intersects $e_1$,~$e_2$ and~$e_3$
(three adjacent edges). Then $e_1$ and $e_2$ have a vertex in common, as
do $e_2$ and $e_3$. Therefore $e_1$ and $e_3$ lie in the plane
determined by $t$ and $e_2$ and thus have vertex in common. But the
three edges cannot form a triangle as each contains a point of~$t$.
\endproof

These conditions provide extra information about non-degenerate
polygons, as well as $n$-secant lines. Condition~(1) shows that a non-degenerate polygon is
indeed a knot. Condition~(2) shows that  no two adjacent edges of~$K$
are collinear and that each component of intersection of an $n$-secant
line with a non-degenerate polygon $K$ is a single point. Condition~(2) also tells us that there are no
multi-vertex trisecant lines. Trisecant lines (and higher order secant
lines) can intersect~$K$ in at most \emph{one} vertex. 

\begin{definition} A generic polygonal knot $K$ has two kinds of trisecants ($n$-secants). A \emph{vertex trisecant ($n$-secant)} includes one
  vertex of $K$ and a \emph{nonvertex trisecant ($n$-secant)} includes no vertices of~$K$.
\end{definition}


Recall that a doubly-ruled surface is a surface with two rulings of
straight lines. Every line of one ruling intersects every line of the
other ruling (maybe at infinity) and each point on the doubly-ruled
surface lies on exactly one line from each ruling. 
It is a well known fact (see for instance \cite{complinegeom} Ch3 or \cite{Otal}) that any three pairwise skew lines generate a
doubly-ruled surface.

\begin{figure}
\centerline{\begin{overpic}[scale=.4]{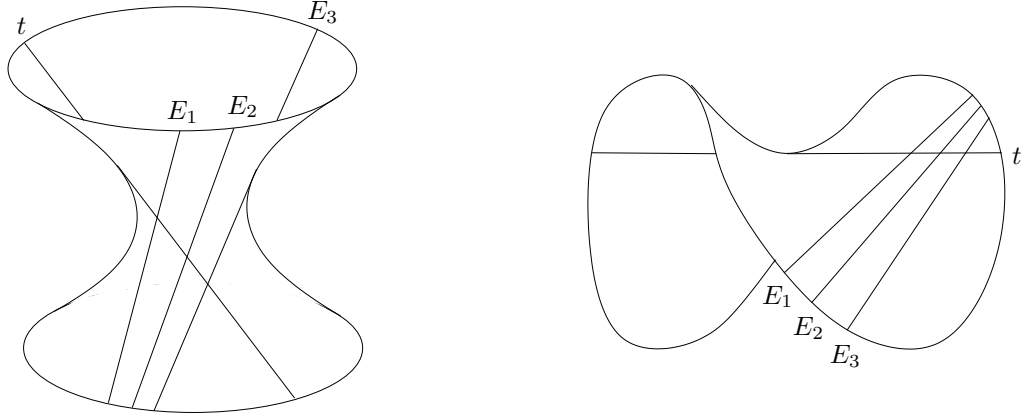}
\put(1,38){$t$}
\put(16,29.5){$E_1$}
\put(22,30){$E_2$}
\put(30,39.5){$E_3$}
\put(75.5,11){$E_1$}
\put(78.6,7.8){$E_2$}
\put(82.3,5){$E_3$}
\put(100.5,25){$t$}
\end{overpic}}
\caption[Doubly-ruled surfaces generated by three pairwise skew
  lines]{Three pairwise skew
  lines determine a doubly-ruled surface either a one sheeted
  hyperboloid (left) or a hyperbolic paraboloid (right). The three skew lines lie
  in one ruling and the lines~$t$ intersecting them lie in the other.} 
\label{fig:qsurface}
\end{figure}

\begin{proposition}\label{prop:surface} A triple of pairwise skew
  lines $E_1$, $E_2$, $E_3$
determines a doubly-ruled surface, either a one-sheeted
hyperboloid or a hyperbolic paraboloid. Moreover, the lines~$E_i$
will lie on one ruling of the surface and any line~$t$ intersecting all the~$E_i$ will lie on the other ruling. See \figr{qsurface}. There is an
open interval or a circle of such lines~$t$. \qed\ 
\end{proposition}

\begin{figure}
\begin{overpic}[scale=.4]{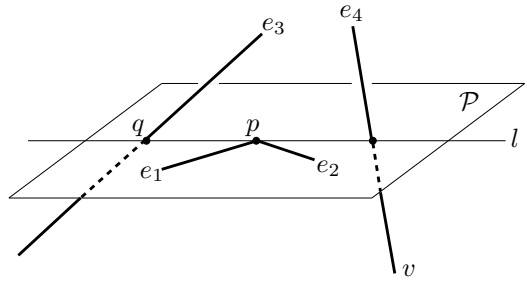}
\put(24,28){$q$}
\put(25.5,19){$e_1$}
\put(46,28){$p$}
\put(59.5,20){$e_2$}
\put(49,47.5){$e_3$}
\put(64,50){$e_4$}
\put(87,32){$\mathcal{P}$}
\put(97,25){$l$}
\put(76,0){$v$}
\end{overpic}
\caption[Trisecant line breaking genericity condition~(G3).]{A
  non-generic trisecant line $l$ lying in the osculating plane of
  vertex $p$, breaking genericity condition (G3) of the
  definition of a generic knot.} \label{fig:nogentri}
\end{figure}

\begin{definition}\label{def:generic} A non-degenerate polygon~$K$ is \emph{generic} if it satisfies the following
  \emph{genericity conditions}:
\item{(G1)} Given any three pairwise skew edges $e_1$, $e_2$, $e_3$ of
  $K$ and the doubly-ruled surface $H$ that they generate, no vertex of $K$
  (except endpoints of $e_1$, $e_2$, $e_3$) may lie in~$H$. 
\item{(G2)} There are no quintisecants (or higher order secants).
\item{(G3)} There are no vertex trisecant lines which lie in the
osculating plane of the vertex.  See \figr{nogentri}.
\end{definition}

\begin{proposition}\label{prop:pfgeneric} The set of all generic
  $n$-gons is open and dense in $\R^{3n}$.
\end{proposition}

\proof The proof uses the same ideas as the proof of \prop{dense} so
is sketched here. In each case, $n$-gons which fail to meet these
conditions lie in an algebraic variety of positive
codimension --- that is non-generic
configurations are closed and nowhere dense. 
These conditions are, in some sense, natural conditions to
consider. If they are broken 
then the lines determined by the edges of the knot lie in degenerate
configurations, as discussed in \cite{transversal}.

(G1) Having fixed $e_1$, $e_2$, $e_3$, the non-generic configurations
are precisely the doubly-ruled surface that they generate: a codimension~1 condition.
(G2) Repeating the discussion in the introduction, we see that secants form a $2$-parameter family, trisecants a $1$-parameter family, quadrisecants a $0$-parameter family and quintisecants a $-1$-parameter family. That is, quintisecants exist only for a codimension~1 set of $n$-gons. (Similarly only a codimension~$k$ set of $n$-gons would have $(k+4)$-secants for~$k>1$.)
(G3) Here, non-generic configurations occur when the fourth edge intersects a line ~$l$ determined by the the two adjacent
edges and the third edge. That is, the plane spanned by~$l$ and one vertex of the fourth edge must also contain the other vertex of the fourth edge.
This is a codimension~1 condition.
\endproof

\section{The Structure of the Set of Trisecants}\label{sect:tri}

The assumption of generic conditions has strong implications for the
structure of the set of trisecants $\mathcal{T}$ in $K^3$ and in
the set of secants $S$. {\bf Unless stated otherwise, all the results of
this section will apply to generic polygonal knotted curves.} 
 
For a generic polygonal knot, no two adjacent edges of $K$ are collinear. Thus~$K^3$ is a union of (closed) cubes $\bigcup_{i,j,k}(e_i\times e_j\times e_k)$. The fattened big diagonal~$\tilde{\Delta}$ in~$K^3$ is a union of cubes along the big diagonal 
 $\Delta := \{abc\in K^3 \,|\, a=b\,\mathrm{or}\, b=c\, \mathrm{or}\, a=c\}$. That is,
 $\tilde{\Delta}= \bigcup_i (e_i\times e_i\times K) \cup\, \bigcup_i  
(e_i\times K\times e_i) \cup\, \bigcup_i (K\times  e_i\times e_i)$. 
As trisecants have three distinct points, no two of which lie on an
edge of~$K$, then~$\mathcal{T}\subset K^3\setm\tilde{\Delta}$. Thinking of $K^3$ as a cubical complex,~$\mathcal{T}$ avoids the 1-skeleton as trisecants do not include more than one vertex of~$K$. Recall that there are two types of trisecants: vertex and nonvertex. Vertex trisecants are just an intersection of~$\mathcal{T}$ with the 2-skeleton of the cubical complex.  Moreover,  as trisecant~$abc$ is an ordered tuple of points then~$abc$ only lies in the cube $e_1\times e_2\times e_3$, where $a\in e_1$, $b\in e_2$ and $c\in e_3$. Trisecant~$abc$ does not lie in $e_2\times e_1\times e_3$, nor in $e_3\times e_2\times e_1$. (However trisecant~$cba$ lies in $e_3\times e_2\times e_1$.) In a similar way,~$K^2$ is a union of squares and the set of secants $S=K^2\setm\tilde{\Delta}$, a topological annulus. 
 
\begin{figure}
\centerline{\begin{overpic}[scale=.4]{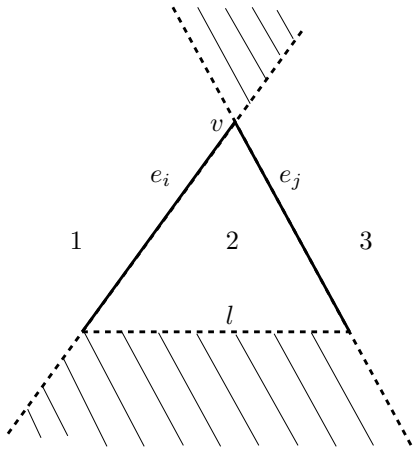}
\put(46.5,71.5){$v$}
\put(33,60){$e_i$}
\put(62,60){$e_j$}
\put(15,45){$1$}
\put(50,45){$2$}
\put(80,45){$3$}
\put(50,28){$l$}
\end{overpic}}
\caption[Regions in the plane spanned by edges $e_i$ and $e_j$ determine the
  type of interval of trisecants.]{Regions in the plane~$\mathcal{P}$
  spanned by edges~$e_i$ and~$e_j$ determine the type of interval of
  trisecants. If the third edge~$e_k$ intersects~$\mathcal{P}$ in the shaded
  regions there are no trisecant lines intersecting $e_i$,~$e_j$,~$e_k$. If $e_k$ intersects~$\mathcal{P}$ in region~2, the set of
  trisecants is homeomorphic to~$[0,1]$. If~$e_k$ intersects~$\mathcal{P}$ in regions~1 or~3, the set of trisecants is
  homeomorphic to~$[0,1)$.}
\label{fig:regions}
\end{figure}

\begin{definition}
Recall that the osculating plane at a vertex~$v$ of a generic polygonal curve is the plane spanned by the two edges incident to~$v$. A \emph{degenerate trisecant} is a triple~$vvp$ or~$pvv$ in $\Delta\subset K^3$, where $p$ is in the osculating plane of vertex~$v$ and the line determined by~$pv$ may not lie in the double wedge bounded by the lines determined by the edges incident to~$v$. (In \figr{regions}, $p$ may only lie in regions 1 or~3.)
\end{definition} 

\begin{definition}
Let $K$ be a generic polygonal curve and $\mathcal{C}=e_1\times e_2\times e_3$ be a cube in~$K^3$. Then we will see that the set $\m{T}\cap\m{C}$ of trisecants in that cube, if nonempty, is homeomorphic to an interval so we call it an \emph{interval of trisecants}. There are no trisecants intersecting three adjacent edges, so either all edges~$e_i$ are pairwise skew and we call it an interval of \emph{skew trisecants} or some~$e_i$ and~$e_j$ are adjacent and we call it an interval of \emph{adjacent trisecants}.
\end{definition}

\begin{lemma}\label{lem:allint}
Let $K$ be a generic polygonal curve. If nonempty, an interval of trisecants~$\m{T}\cap\m{C}$ in a cube~$\mathcal{C}=e_1\times e_2\times e_3$ in~$K^3$ is homeomorphic to a closed or half-open interval. The endpoint or endpoints lie in (the interior of) distinct faces of~$\m{C}$, and the interior of the interval lies in the interior of~$\m{C}$. The half-open interval only occurs when~$e_2$ is adjacent to either~$e_1$ or~$e_3$; here, the open end of the interval approaches a degenerate trisecant in the (interior of) an edge of~$\m{C}$. Along an interval of skew trisecants, each of the three points~$p_i\in e_i$ moves monotonically and smoothly along the edge~$e_i$. Along an interval of adjacent trisecants, if $e_i$ is the nonadjacent edge, then~$p_i$ is constant, while the other~$p_j$ move smoothly and monotonically along~$e_j$.
\end{lemma}

\begin{figure}
\centerline{\begin{overpic}[scale=.5]{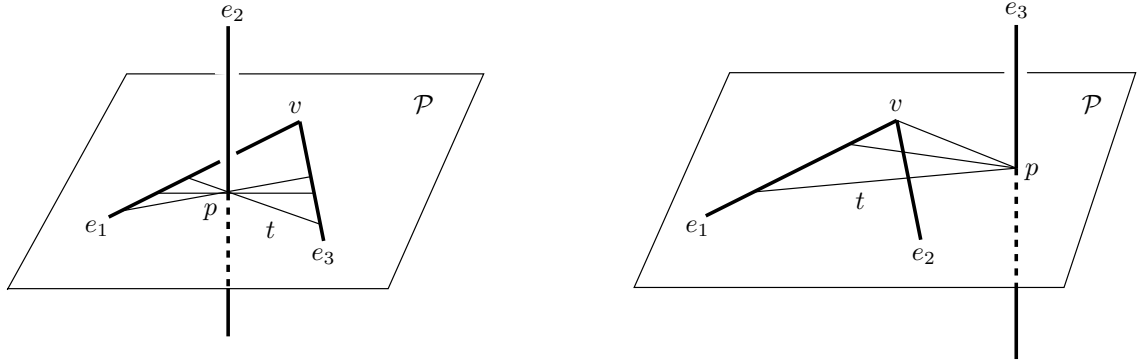}
\put(7,11.5){$e_1$}
\put(27,9){$e_3$}
\put(19,30.5){$e_2$}
\put(25,22){$v$}
\put(23,11){$t$}
\put(17.5,13){$p$}
\put(36,22){$\mathcal{P}$}
\put(60,11.5){$e_1$}
\put(80,9){$e_2$}
\put(88.2,30.8){$e_3$}
\put(78,22.2){$v$}
\put(75,13.5){$t$}
\put(90,16.5){$p$}
\put(95,22){$\mathcal{P}$}
\end{overpic}}
\caption[Intervals of adjacent trisecants.]{Intervals
  of adjacent trisecants. The left picture shows an
  interval homeomorphic to a closed interval, the right an
  interval homeomorphic to a half-open interval.} 
\label{fig:casetwo}
\end{figure}

\proof 
We first consider where $\m{T}\cap\m{C}$ is an interval of adjacent trisecants.
Let~$\mathcal{P}$ be the plane spanned by two adjacent edges~$e_i$ and~$e_j$ and let~$p$ be the unique point of intersection of the third edge~$e_k$ with~$\mathcal{P}$. (Non-degeneracy implies that~$p$ is in the interior of~$e_k$.)
\figr{regions} shows~$\mathcal{P}$ separated into regions divided by
lines determined by~$e_i$ and~$e_j$ and the line~$l$ between the non-intersecting
vertices of~$e_i$ and~$e_j$. For trisecants to occur then point~$p$
must be in one of regions 1, 2 or~3 (and not the shaded regions).  There is a
1-parameter family of such trisecants. If~$p$ is in region~2, then~$e_1$ is adjacent to~$e_3$ and trisecants are of the form~$apb$, where $a\in e_1$ and $b\in e_3$. Clearly, the interval of trisecants is homeomorphic to a closed interval as in \figr{casetwo}~(left). If~$p$ is in either region~1 or~3, then~$e_2$ is adjacent to either~$e_1$ or~$e_3$ and trisecants are of the form~$abp$ where $a\in e_1$ and $b\in e_2$ or $pab$ where $a\in e_2$ and $b\in e_3$. \figr{casetwo}~(right) shows the first possibility. In either case, the interval of
trisecants will be homeomorphic to a half-open interval and ends on the degenerate trisecant $vvp$ (or $pvv$) where
$v$ is the common vertex of $e_1$ and $e_2$ (or $e_2$ and $e_3$).  This point is not in the set of trisecants, $vvp$ (or $pvv$) $\in\Delta\subset\tilde{\Delta}$.

\begin{figure}
\centerline{\begin{overpic}[scale=.4]{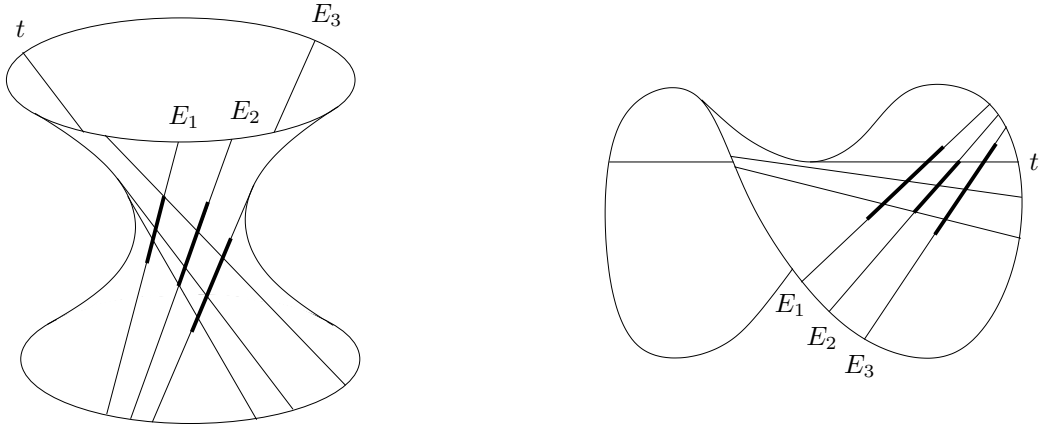}
\put(1,38){$t$}
\put(16,29.5){$E_1$}
\put(22,30){$E_2$}
\put(30,39.5){$E_3$}
\put(75.5,11){$E_1$}
\put(78.6,7.8){$E_2$}
\put(82.3,5){$E_3$}
\put(100.5,25){$t$}
\end{overpic}}
\caption[Closed interval of trisecant lines on doubly-ruled
surfaces.]{Two doubly-ruled surfaces, each with a closed interval of trisecant lines intersecting edges~$e_1$,~$e_2$ and~$e_3$.}
\label{fig:interval}
\end{figure}

We now consider where $\m{T}\cap\m{C}$ is an interval of skew trisecants.
By \prop{surface}, the triple of pairwise skew edges
determines a doubly-ruled surface, 
either a one-sheeted hyperboloid or hyperbolic paraboloid. The edges
$e_i$~$(i=1,2,3)$ all lie in one ruling. A trisecant line is a line in
the other ruling which intersects~$e_1$,~$e_2$ and~$e_3$. Let~$E_i$ be
the line determined by~$e_i$. From \prop{surface} we know there is an
circle of lines intersecting~$E_1$,~$E_2$ and~$E_3$. Let~$I_1$ be the closed sub-interval of these lines which intersects the
edge~$e_1$ and the lines~$E_2$ and~$E_3$. As~$e_1$ is an edge
(closed),~$I_1$ is homeomorphic to a closed interval.  Let~$I_2$ and~$I_3$ be defined in a similar manner, they too are homeomorphic to a closed interval. 
Thus, the set of trisecants though~$e_1$,~$e_2$ and~$e_3$ is homeomorphic to the
intersection of~$I_1$,~$I_2$ and~$I_3$ (three closed intervals in an
open interval), which is either a closed interval or a point, see
\figr{interval}. Note that this interval of trisecants must start and
end with a vertex trisecant corresponding to one of the
endpoints of the~$I_i$. The interval cannot have zero length (be a point) nor start with a line intersecting two (or more) endpoints of~$I_i$, as either of these cases corresponds to a trisecant
with two (or more) vertices of the knot --- contradicting
non-degeneracy.  Hence, if it is nonempty, the set of
trisecants through~$e_1$,~$e_2$ and~$e_3$ is homeomorphic to a closed
interval of non-zero length which starts and ends with vertex
trisecants. 

We may say more about the structure of~$\mathcal{T}\cap \m{C}$ for both intervals of adjacent and skew trisecants. In all cases the interior of the interval of trisecants lies in the interior of~$\m{C}$ and as the endpoints correspond to vertex trisecants, they lie on the interior of a face of~$\m{C}$. Parameterize~$K$ with respect to arclength; along an interval of skew trisecants, each of the three points~$p_i\in e_i$ moves monotonically and smoothly along the edge~$e_i$.  Along an interval of adjacent trisecants, on the nonadjacent edge~$p$ is constant while on the adjacent edges, the points move monotonically and smoothly along their corresponding edges.
 (An explicit calculation for both cases may be found in \cite{Den} in the Appendix.)  
Thus~$\m{T}\cap\m{C}$ is smoothly embedded in~$\m{C}$.
\endproof

\begin{proposition}\label{prop:embedded}
Let $K$ be a generic polygonal curve. If nonempty, ~$\overline{\mathcal{T}}$ is a compact 1-manifold with boundary,
embedded in~$K^3$ in a piecewise smooth way with $\mathcal{T}\subset
K^3 \setm \tilde{\Delta}$ and
$\partial{\mathcal{T}}\subset{\Delta}$. Moreover~$\bdry{\m{T}}$ is exactly the set of degenerate trisecants.
\end{proposition}

\proof The definition of a trisecant shows $\mathcal{T}\subset
K^3 \setm \tilde{\Delta}$. Consider a face~$\m{F}$ between two cubes~$\m{C}$ and~$\m{C}'$ in~$K^3$. Suppose the interval of trisecants in~$\m{C}$ ends at a particular vertex trisecant~$t$ in~$\m{F}$. There are three possibilities for the interval of trisecants in~$\m{C}'$: there is an interval of trisecants starting at~$t$, there is an interval of trisecants starting at another vertex trisecant~$t'$ in~$\m{F}$ or there are no intervals of trisecants which start on~$\m{F}$.
In the last two cases, the intersection $\m{T}\cap\m{C}'$ is either a single point or an interval of trisecants and a single point. \lem{allint} rules out these cases. Thus in~$\mathcal{C}'$, there must be an interval of trisecants starting at~$t$. Hence~$\mathcal{T}$ is a 1-manifold embedded
in $K^3$ in a piecewise smooth way.

To understand the claim about $\overline{\mathcal{T}}$, we must
understand~$\bdry\mathcal{T}$. Recall that the half open intervals of
trisecants occur when an interval of adjacent trisecants ends on a
degenerate trisecant~$vvp$ (or~$pvv$) where~$v$ is a vertex of~$K$. Thus the limit point (what is added in~$\overline{\mathcal{T}}$) is
a degenerate trisecant and lies in~$\Delta$. Also note that no other interval of trisecants in~$K^3$ can end at the degenerate trisecant~$vvp$ (or~$pvv$). Degenerate trisecant~$vvp$ can be a common degenerate vertex trisecant to only two cubes in~$K^3$. Without loss of generality denote these cubes $e_1\times e_2\times e_3$  and $e_2\times e_1\times e_3$ where~$v$ is the common vertex of~$e_1$ and~$e_2$. But trisecants are ordered tuples, so at most one of these has nonempty intersection with~$\m{T}$. Similarly for~$pvv$. Hence~$\overline{\mathcal{T}}$ is a compact 
1-manifold with boundary embedded in a piecewise smooth way in~$K^3$.
\endproof


We now wish to understand the structure of $T=\pi_{12}(\mathcal{T})$ the projection of
the set of trisecants to the set of secants $S=K^2\setm\tilde{\Delta}$. Recall that a secant~$ab\in S$ is only in~$T$ if it is the first two points of trisecant~$abc$.

We first put a metric on~$S$ (and hence a topology) in the following way.
The knot~$K$ has an orientation, is parameterized with respect to
arclength and~$d(a,a')$ is the shorter arclength between~$a$ and~$a'$ along~$K$. We give~$K^2$
the product metric such that $d((a,b),(a',b'))= d(a,a')+d(b,b')$. Now consider~$S$ as a metric space, not with the subset metric, but with the path metric $d_S((a,b),(a',b'))=\inf_{\gamma}(len(\gamma))$, where~$\gamma$ is any
path in~$S$ from~$(a,b)$ to~$(a',b')$ and~$len(\gamma)$ the
length of~$\gamma$ (using the product metric from $K^2$). Let the completion of~$S$ in
this metric be  
$$\overline{S} =(K^2\setm\tilde{\Delta}) \cup \tilde{\Delta}_+ \cup
\tilde{\Delta}_-\, .$$

Here $\tilde{\Delta}_-$ denotes the lower edge of $S$, that is
approaching $\tilde{\Delta}$ in $K^2$ from above. That is $ab\in K^2$ and
$b>a$: $a$ and $b$ are close together with~$b$ just after~$a$ in the order of
$K$. Similarly $\tilde{\Delta}_+$ denotes the upper edge of $S$. Let~$\Delta_-$,~$\Delta_+$ denote the lower
and upper big diagonals of~$S$ respectively. The topology of~$S$ is
imposed by  the metric. Hence the upper and lower edges of the annulus~$S$ are far from each other. The shortest path between them lies \emph{in}~$S$ and by definition can \emph{not} cross~$\Delta$.

Let $\overline{T}$ be the projection of $\overline{\mathcal{T}}$ in
$\overline{S}$, $\overline{T}:=\pi_{12}(\overline{\mathcal{T}})$.
Now $\pi_{12}$ maps to $K^2$ and $\overline{S}\nsubseteq K^2$, thus to
understand $\overline{T}$ we
must describe where $\pi_{12}(\bdry T)$ lies in $\overline{S}$.

\begin{figure}
\centerline{\begin{overpic}[scale=.4]{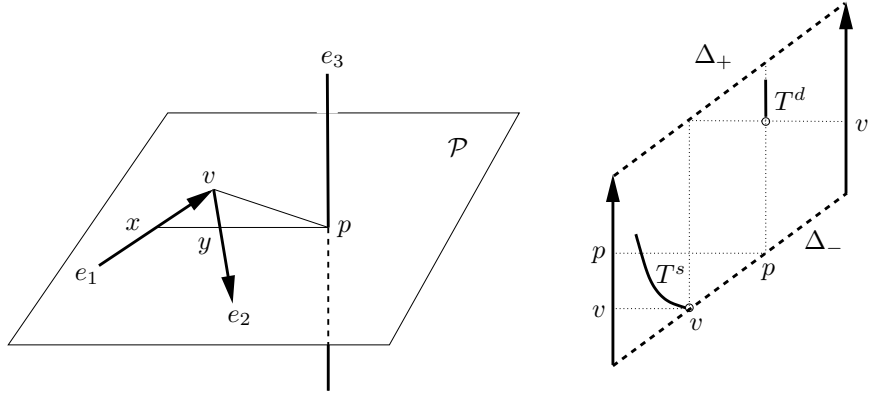}
\put(8,13.4){$e_1$}
\put(26,8.5){$e_2$}
\put(37,39){$e_3$}
\put(52,28){$\mathcal{P}$}
\put(39,19){$p$}
\put(14,19.3){$x$}
\put(22.5,17){$y$}
\put(23,24.8){$v$}
\put(76.5,13){$T^s$}
\put(69,16){$p$}
\put(69,9){$v$}
\put(80.5,7.5){$v$}
\put(89,14){$p$}
\put(94,17){$\Delta_-$}
\put(81,39){$\Delta_+$}
\put(90.5,33.5){$T^d$}
\put(100,31){$v$}
\end{overpic}}
\caption[The interval of adjacent trisecants ends on a degenerate
trisecant.]{The left picture shows an interval of adjacent trisecants
  ending on degenerate trisecant $vvp$ or $pvv$. The right picture shows
  the corresponding intervals of trisecants in $S$. The intervals of
  $T^s$ and $T^d$ correspond to trisecants with 
  linear ordering $e_1e_2e_3$ and $e_3e_2e_1$ respectively.}
\label{fig:position} 
\end{figure}

A half-open interval of trisecants occurs when an interval of adjacent
trisecants ends on a degenerate trisecant. (See
Figures~\ref{fig:casetwo} and~\ref{fig:position}.) Let the points of
the adjacent trisecants lie on edges~$e_1$,~$e_2$ and~$e_3$ and let~$e_1$ and~$e_2$ be adjacent with common
vertex~$v$. Let~$e_3$ intersect the plane spanned by~$e_1$ and~$e_2$
in the point~$p$. Suppose the
points of the trisecants lie on~$e_1e_2e_3$ in that order, then the
third point~$p$ is fixed. Thus the interval of trisecants ends at
degenerate trisecant~$vvp\in\Delta\subset K^3$ and in $\overline{S}$,
the interval ends at either~$\Delta_+$ or~$\Delta_-$. Now suppose the points of
the trisecant lie on~$e_3e_2e_1$ in that order, then the first point~$p$ is fixed. Thus the interval of trisecants ends at degenerate trisecant
$pvv\in\Delta\subset K^3$ and in $\overline{S}$, it ends at~$pv$ in the
interior of~$S$.

In fact we can say more.  First consider the case where knot is
oriented from~$e_1$ to~$e_2$ as illustrated in \figr{position} (left). Trisecants whose points lie on edges~$e_1e_2e_3$ in that order are trisecants of same ordering~$\m{T}^s$. Moreover, as the interval of trisecants~$abp$ ends at
degenerate trisecant~$vvp$, the first point of the trisecant~($a\in
e_1$) moves with the parameterization of~$K$ and the
second~($b\in e_2$) against. Thus in~$K^2$, $b>a$ and the interval of trisecants approaches the diagonal from above. Hence this part of $T^s$ has negative slope and ends at~$vv\in\Delta_-$. In
this situation, $\pi_{12}(\bdry\mathcal{T})\subset \Delta_-$. Trisecants whose points lie on~$e_3e_2e_1$ in that 
order are trisecants of different ordering~$\m{T}^d$. The first point of the trisecant is fixed and
the second moves against the parametrization. Hence this interval of
trisecants is a vertical line, decreasing and ending at the point
$pv\in S$ and $\pi_{12}(\bdry\mathcal{T})\subset S$. \figr{position}
(right) shows the two corresponding intervals of trisecants in~$S$.

Now repeat these arguments but with the orientation of $K$
reversed --- from~$e_2$ to~$e_1$. Trisecants whose points lie on edges~$e_1e_2e_3$ in that order are now in~$T^d$ and the interval of
trisecants has negative slope and ends on $vv\in\Delta_+$ and  $\pi_{12}(\bdry\mathcal{T})\subset \Delta_+$. Trisecants whose points
lie on edges $e_3e_2e_1$ in that order are in $T^s$ and the interval of trisecants is a vertical line increasing and ending at $pv\in S$ and $\pi_{12}(\bdry\mathcal{T})\subset S$.

We have completely described the behavior of $\bdry T:=\pi_{12}(\bdry
\mathcal{T})$ in~$\overline{S}$. Either~$\bdry T$ lies in $S=K^2\setm \tilde{\Delta}$ as the 
limit point of a vertical interval of trisecants whose first
point is fixed, or~$\bdry T$ lies on~$\Delta_+$ or~$\Delta_-$ as the limit point of a negatively sloped interval of
trisecants whose third point is fixed. 

\begin{definition}
Recall $\overline{T}=\pi_{12}(\overline{\mathcal{T}})$ and $\bdry T:=
\pi_{12}(\bdry \mathcal{T})$. In a similar way, 
define $\overline{T}^s=\pi_{12}(\overline{\mathcal{T}}^s)$,
$\overline{T}^d=\pi_{12}(\overline{\mathcal{T}}^d)$, $\bdry T^s=
\pi_{12}(\bdry \mathcal{T}^s)$ and $\bdry T^d= \pi_{12}(\bdry \mathcal{T}^d)$.
\end{definition}

In \lem{intint} we will show that 
if $\overline{T}^s \cap \overline{T}^d\neq\emptyset$, then $T^s\cap T^d
\neq\emptyset$. Hence by \lem{alt} it will be sufficient to show that 
 $\overline{T}^s \cap \overline{T}^d\neq\emptyset$ in $\overline{S}$
in order to show that an alternating quadrisecant exists.

\begin{lemma}\label{lem:immersed} Let $K$ be a generic polygonal
  curve. The projection $\pi_{12}$ is a piecewise smooth
  immersion of $\mathcal{T}$ into $S$.
\end{lemma}
 
\proof Let $\m{C}$ be a cube in $K^3$. From \lem{allint}, $\mathcal{T}\cap\m{C}$ is both smooth and monotonic in at least two variables (corresponding to edges of~$\m{C}$). Thus the tangent vector is never vertical. The map $\pi_{12}(\mathcal{T})$ is just an orthogonal projection hence the projection~$\pi_{12}$ of $\mathcal{T}\cap\m{C}$ in~$S$ is a smooth immersion.  If nonempty, $\mathcal{T}$ is an embedded 1-manifold in $K^3$. Thus~$\pi_{12}$ is a
piecewise smooth immersion of~$\mathcal{T}$ in~$S$.
\endproof

\begin{definition}  Let $K$ be a generic polygonal curve and let $\m{C}$ be a cube in~$K^3$. The set  $\pi_{12}( \mathcal{T}\cap\m{C})$, if nonempty, is homeomorphic to an interval. Thus we (again) call it an \emph{interval of trisecants} or an \emph{interval of~$T$ in~$S$}.
\end{definition}

The following two lemmas show that in $S$, $T$ intersects itself
only in double-points. This means that quadrisecants (if they exist)
are isolated. 

\begin{lemma}\label{lem:onequad} Let $K$ be a generic polygonal
  curve. In~$S$, the intersection of two intervals of~$T$
   is either be empty, one point, or two distinct points. 
\end{lemma}

\proof We assume that $T$ is nonempty or the result is trivial. First assume that the two intervals of~$T$ in~$S$ do not share a common vertex. The two intervals have
a point in common if and only if two trisecants have the same
first and second points. That is, only if they have the same first and
second edges. However, the third edges differ, or else there would not
be two intervals of~$T$ intersecting in~$S$. There are three cases to consider. Let
the first two edges be denoted~$e_1$ and~$e_2$ and let the differing
third edges be denoted by~$e_3$ and~$e_3'$.

\begin{figure}
\centerline{\begin{overpic}[scale=.4]{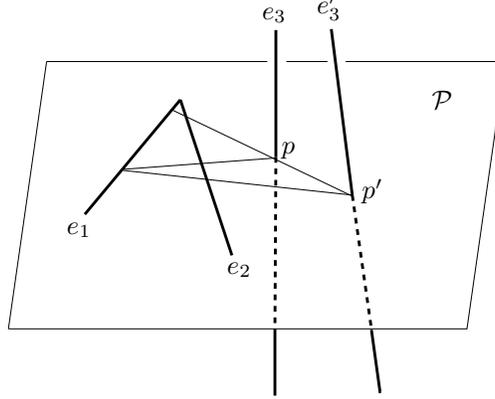}
\put(12,33){$e_1$}
\put(44,25){$e_2$}
\put(51,76){$e_3$}
\put(62,77){$e_3'$}
\put(85,58){$\mathcal{P}$}
\put(55,49){$p$}
\put(71,40){$p'$}
\end{overpic}}
\caption{Trisecants with same two adjacent first edges can have at
  most one line in common.}\label{fig:triplane}
\end{figure}

{\bf Case 1:} Both intervals are intervals of adjacent trisecants. There is
  only one possible case, the first two edges are adjacent. (If $e_1$
  and~$e_3$ are adjacent and~$e_1$ and~$e_3'$ are adjacent, then a
  quadrisecant whose points lie on $e_1e_2e_3e_3'$ includes a
  trisecant whose points lie on three adjacent
  edges $e_1e_3e_3'$, contradicting non-degeneracy.) Both
  intervals of trisecants lie in the plane~$\mathcal{P}$
  spanned by~$e_1$ and~$e_2$. Let $e_3$ intersect~$\mathcal{P}$ in the
  point~$p$ and~$e_3'$ intersect~$\mathcal{P}$ in the
  point~$p'$. \figr{triplane} illustrates this. The line through~$pp'$ is the only possible common trisecant line. 

{\bf Case 2:} One interval is an interval of adjacent
  trisecants and the other an interval of skew trisecants. Let the points of the
  skew trisecant lie on edges $e_1$, $e_2$, $e_3$; these determine a
  doubly-ruled surface. Let the
  points of the adjacent trisecants lie on edges~$e_1$,~$e_2$ and~$e_3'$. As~$e_1$ and~$e_2$ are pairwise skew, either~$e_1$
  and~$e_3'$ are adjacent, or~$e_2$ and~$e_3'$ are adjacent. In either case, 
  the plane spanned by the two adjacent edges intersects the
  doubly-ruled surface in some quadratic curve. However, as one of
  the edges lies in the doubly-ruled surface, the curve
  degenerates to a pair of straight lines, one of which includes edge~$e_1$ (or~$e_2$). The other straight line is the only possible
  common trisecant line and only if it is both a trisecant line on the
  doubly-ruled surface and the plane. 
  
{\bf Case 3:} Both intervals are intervals of skew
  trisecants. Let~$H$ be the 
  doubly-ruled surface generated by the lines~$E_i$ determined by the
  polygonal edges~$e_i$. Suppose~$e_{3}'$ intersects~$H$,
  it does so in one or two points. The trisecant
  lines(s) through these point(s) will 
  only correspond to a quadrisecant line if it is a trisecant through the
  edges $e_1$, $e_2$, $e_3$, \emph{and} through the edges $e_1$,
  $e_2$, $e_3'$. In fact this is all that can happen. The edge~$e_3'$ can
  never be contained in~$H$. If it were
  then its vertices lie in~$H$ which
  contradicts genericity condition (G1). 
Thus in all three cases the intersection of two
intervals of~$T$ is either empty, or one or two points. 

Finally assume the two intervals of~$T$ have a common vertex. This corresponds to a vertex trisecant which joins two intervals of trisecants.  In~$K^3$, these intervals of trisecants lie in adjacent cubes $\m{C} = e_1\times e_2 \times e_3$ and $\m{C'}=e' _1\times e'_2 \times e'_3$.  As the cubes are adjacent, two of the edges are the same and the third differs. 
If $e_1\neq e'_1$ (or $e_2\neq e'_2$) then the first (or second) points of
the intervals of trisecants lie on different edges. Thus the only possible common
point in~$S$ is the projection of the common vertex trisecant.
If the first two edges are the same and the third differs, then the intervals of trisecants may intersect in~$S$. There are three cases to consider as above. The details are omitted as they are similar to the previous arguments. 
\endproof

\begin{lemma}\label{lem:doublepoints} Let $K$ be a generic polygonal curve. No more than two points of~$\mathcal{T}$ can have the same image under~$\pi_{12}$. 
\end{lemma}

\begin{figure}
\centerline{\begin{overpic}[scale=.4]{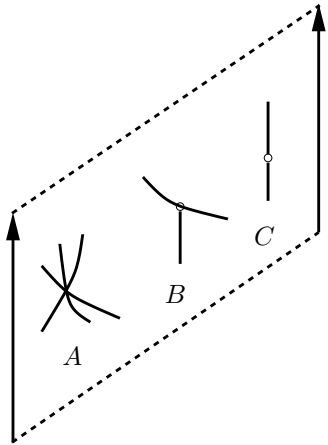}
\put(13,18){$A$}
\put(36,32){$B$}
\put(56,45){$C$}
\end{overpic}}
\caption{Intersections of $T$ in $S$ forbidden for a generic polygonal
  knot.}\label{fig:triint} 
\end{figure}

\proof  If there were three (or more) points of~$\mathcal{T}$ with the same
image under~$\pi_{12}$, this would 
correspond to a quintisecant (or higher order secant) in contradiction
to genericity condition~(G2). (See also
\figr{triint}A.) 
\endproof

The structure of~$T$ described is sufficient to prove that
nontrivial generic polygonal knotted curves have an alternating
quadrisecant. However, with the addition of extra genericity conditions, it is possible to provide even more details about the structure of~$T$. These genericity conditions and details are found in~\cite{Den}. 

Let $\pi_{ij}(K^3)\rightarrow K^2$ be projection
onto the $i$th and $j$th coordinates (where $i<j$ and $i,j=1,2,3$.)
We have extensively studied the structure of~$\pi_{12}(\mathcal{T})$
for generic polygonal curves and have shown that, if nonempty, it is the image of a
piecewise smooth immersion and which intersects itself
at double-points. We can also say the same thing about
$\pi_{13}(\mathcal{T})$ and $\pi_{23}(\mathcal{T})$. \lem{immersed}
carries over immediately. Similarly \lem{onequad} and \lem{doublepoints} can be proved
for $\pi_{13}(\mathcal{T})$ and $\pi_{23}(\mathcal{T})$ with only
minor alterations. Thus the following proposition holds:

\begin{proposition}\label{prop:allproj} Let $K$ be a generic polygonal curve.
In $K^2$, the projection~$\pi_{ij}$ ($i<j$ and
$i,j=1,2,3$) is a piecewise smooth immersion of~$\mathcal{T}$ into~$S$
and $T=\pi_{ij}(\mathcal{T})$ intersects itself at
double points. \qed
\end{proposition}

For the most part we are only interested in $T=\pi_{12}(\mathcal{T})$
in~$S$. Recall that $\pi_{12}(\mathcal{T})$ allows us to capture the
orderings of  coincident trisecants which in turn enables us to find
alternating quadrisecants. In \lem{intint} we will show that if
$\overline{T}^s\cap 
\overline{T}^d\neq\emptyset$ in~$\overline{S}$, then $T^s\cap
T^d\neq\emptyset$ in~$S$. Thus it is important to
understand exactly where $\overline{T}^s\cap
\overline{T}^d\neq\emptyset$ but~$T^s$ and~$T^d$ do not intersect. Recall
$\overline{T}=\pi_{12}(\overline{\mathcal{T}}) =\pi_{12}(\mathcal{T}
\cup \bdry \mathcal{T})$. After \prop{embedded} we showed that
$\bdry T=\pi_{12}(\bdry \mathcal{T})$ lies on~$\Delta_+$ or~$\Delta_-$ as the
limit point of an interval of trisecants whose first two points both
end on a vertex of~$K$ and whose third point is fixed. If the second
and third points of the trisecants both end on a vertex of~$K$ and the
first is fixed, then~$\bdry{T}$ lies in~$S$.

The following lemma shows that $\overline{T}^s$ stays away
from~$\tilde{\Delta}_+$ and~$\overline{T}^d$ stays away from~$\tilde{\Delta}_-$. Hence 
$\overline{T}^s$ and $\overline{T}^d$ cannot intersect on
$\Delta_+$ or $\Delta_-$, but only in~$S$. 

\begin{figure}
\centerline{\begin{overpic}[scale=.4]{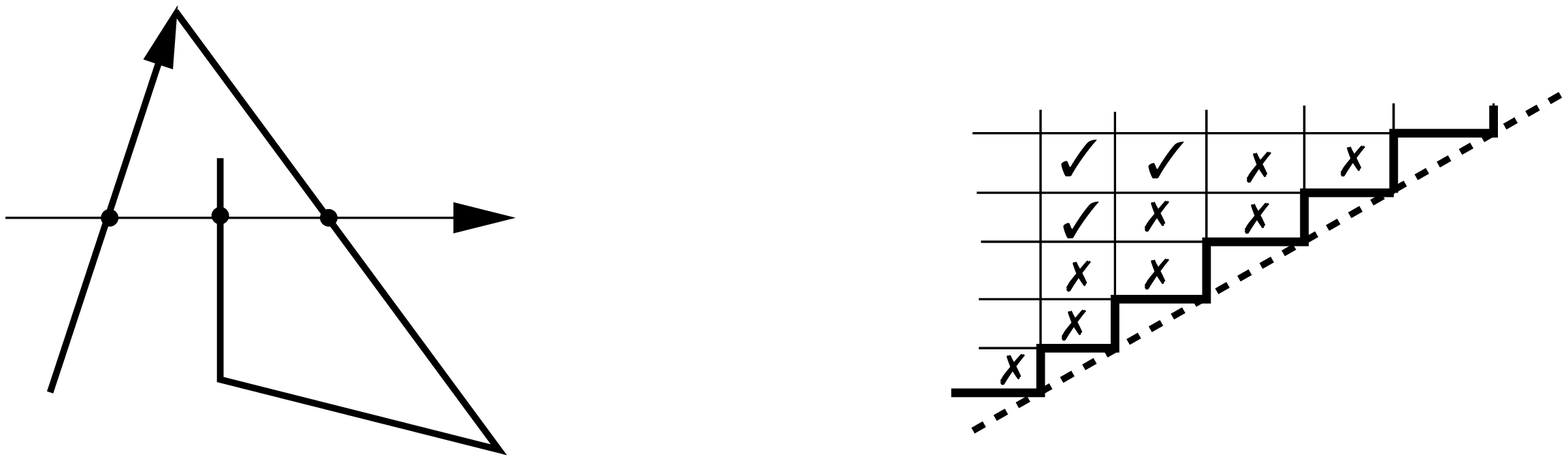}
\put(7,12.5){$a$}
\put(14.5,12){$b$}
\put(23.5,13){$c$}
\put(75,25){$S$}
\put(80,8){$\Delta_-$}
\put(56,3){$\tilde{\Delta}_-$}
\end{overpic}}
\caption[The first two points of trisecants of different order must be at
  least two edge lengths apart.]{The first two points of trisecants of
  different order must be at least two edge lengths apart. This
  information is shown on part of the set of secants $S$ on
  the right. The ticked boxes indicate where trisecants of different
  order may be.}\label{fig:dist} 
\end{figure}

\begin{lemma}\label{lem:dist}
Let $K$ be a generic polygonal knotted curve. In $S$,
  $d_S(\overline{T}^d,\tilde{\Delta}_-)\geq h$ and
  $d_S(\overline{T}^s,\tilde{\Delta}_+)\geq h$, where $h$ is the minimum edge
  length of $K$.
\end{lemma}

\proof  As secant $ab$ approaches $\tilde{\Delta}_-$ it is the same as approaching the diagonal in $K^2$ from above. Here,~$a$ and~$b$ are close together and $b$ is just  after~$a$ in the order of the knot. Now consider secant $ab\in T^d$. The order along the knot is~$acb$ (where~$c$ is the third
point of the trisecant.) Thus~$a$ and~$b$ are
separated by at least one edge as~$c$ must come in between them and by
definition, distinct points of trisecants lie on distinct edges. In
fact~$a$ and~$b$ must be separated by at least two edges. If they were
just separated by one edge, then the points of the trisecant~$abc$ lie
on three consecutive edges contradicting
non-degeneracy. \figr{dist}~(left) illustrates this and 
\figr{dist}~(right) shows this information recorded on~$S$. Trisecants
of different order may be in boxes marked with a tick, not a cross. They are at least one edge away
from $\tilde{\Delta}_-$ (indicated by a bold line). Hence
$d_S(T^d,\tilde{\Delta}_-)\geq h$, where $h$ is the minimum edge length of $K$. A similar argument shows~$d_S(T^s,\tilde{\Delta}_+)\geq h$. 
\endproof

\begin{lemma}\label{lem:intint} Let $K$ be a generic polygonal
  curve. If $\overline{T}^s\cap \overline{T}^d\neq\emptyset$ in
  $\overline{S}$, then $T^s\cap T^d\neq\emptyset$ in $S$.
\end{lemma}

\proof \lem{dist} shows $\overline{T}^s$ and $\overline{T}^d$ can
only intersect in the interior of~$S$. Suppose
$t\in\overline{\mathcal{T}}^s$ and $t'\in\overline{\mathcal{T}}^d$
with $\pi_{12}(t)=\pi_{12}(t')\in S$. We want to show 
$t\in{\mathcal{T}}^s$ and $t'\in{\mathcal{T}}^d$.
Assume, by way of contradiction that $t\in\bdry\mathcal{T}^s$ and
consider two cases: $t'\in\mathcal{T}^d$ and $t'\in\bdry\mathcal{T}^d$.

{\bf Case 1:} Assume that $t\in\bdry\mathcal{T}^s$ and
$t'\in\mathcal{T}^d$ and $\pi_{12}(t)=\pi_{12}(t')\in
S$. \figr{triint}B illustrates this case. 

Let $pv=\pi_{12}(t)=\pi_{12}(t')$. As $t\in\bdry\mathcal{T}^s$, in~$S$
there is a vertical interval of trisecants which increases to and
ends at secant $pv$.  Thus $pv$ is the projection of degenerate trisecant $pvv\in\m{T}$.  Let~$p$ belong to edge~$e_1$ and vertex~$v$ belong to edges~$e_2$ and~$e_2'$.
The vertical open interval of trisecants corresponds to
trisecants with first point~$p$ and second and third points on edges~$e_2$ and~$e_2'$. But
$pv=\pi_{12}(t')$, where $t'\in\mathcal{T}^d$. Thus~$pv$ is also the
first two points of a trisecant of different order~$pva$, where~$a$
lies on edge~$e_3$. Secant~$pv$, and hence trisecant~$pva$ lies in the
plane spanned by~$e_2$ and~$e_2'$. Thus trisecant~$pva$ lies in
the osculating plane of vertex~$v$, contradicting genericity
condition~(G3).

{\bf Case 2:} Assume that $t\in\bdry\mathcal{T}^s$ and
$t'\in\bdry\mathcal{T}^d$ and $\pi_{12}(t)=\pi_{12}(t')\in
S$. \figr{triint}C illustrates this case.

Let $pv=\pi_{12}(t)=\pi_{12}(t')$. As $t\in\bdry\mathcal{T}^s$ and
$t'\in\bdry\mathcal{T}^d$, $pv$ belongs to two vertical intervals
of trisecants which both end at secant~$pv$ in~$S$. The
trisecants of same and different order both have~$p$ as a starting
point and both have points which lie on edges which have~$v$ as a
common vertex. But there can only be two edges with~$v$ as a common
vertex. This implies~$t$ and~$t'$ have the same first point and both
have second and third points on the same two edges. This implies that
the two intervals of trisecants are in fact the same (and hence have the same
order). This contradicts the fact~$t$ and~$t'$ are trisecants of same
and different order respectively. 
\endproof


In \lem{alt} we showed that an alternating quadrisecant exists when
$T^s$ and $T^d$ have a common point in $S$. In \sect{genpoly} and
\sect{tri} we examined the structure of the set of trisecants for
generic polygonal knots. Through  \lem{intint} we see that it is
sufficient to prove $\overline{T}^s$ and $\overline{T}^d$ have a
common point in~${S}$ in
order to show that an alternating quadrisecant exists. In fact, we
will prove a stronger result than this in \sect{quadgenpoly}. 

\section{Essential Quadrisecants}\label{sect:equad}

Eventually we want to prove that all nontrivial tame
knots have an alternating quadrisecant. To do this we will prove that
nontrivial generic polygonal knots have an \emph{essential} alternating
quadrisecant and then use a limit argument. We will define essential
shortly --- this is the topological notion required so the limit
argument works. 

In \sect{genpoly} we showed that the set of generic polygonal knots is
open and dense in the set of polygonal knots. For any knotted curve
$K$, there is a sequence of generic polygonal knots
$\{K_i\}_{i=1}^{\infty}$ which 
converge to $K$. In \lem{gpconverge} we show how this
sequence of $K_i$ may be chosen so that $K_i$ is ambient isotopic to~$K$ and converges in~$C^0$ to~$K$.

Thus we have $K_i\rightarrow K$. \thm{equad} shows that each generic
polygonal knot~$K_i$ has an essential alternating quadrisecant
$a_ib_ic_id_i$. As each $K_i$ is contained in some ball in $\R^3$,
eventually the entire sequence $\{K_i\}_{i=1}^\infty$ lies in a common
compact subset of 
$\R^3$. By taking a subsequence if necessary, there is a
sequence of essential alternating quadrisecants which converge to some
quadrisecant $abcd\in K^4$. We may also assume that the order of
$a_ib_ic_id_i$ along $K_i$ is $a_ic_ib_id_i$.

\begin{lemma}\label{lem:converge} Given the sequence of generic
  polygonal knots with converging alternating quadrisecants as described
  above. If $b$ does not lie on the same subarc of~$K$ as $c$, then no two
  limit points lie on the same straight subarc of~$K$ and so $abcd$ is
  a quadrisecant. 
\end{lemma}

\proof Given two points $p,q\in K$, we denote $p\sim q$ to mean $p$
and $q$ lie on the same common straight subarc of $K$. Quadrisecant
points lie on the same straight subarc of $K$ in the
limit if and only if they are next to each other in the ordering of
{\it both} the quadrisecant line and the knot. Thus $b\nsim c$ implies
$a\nsim b$ as the order along the knot is $a_ic_ib_i$. Similarly
$b\nsim c$ implies $c\nsim d$ as the knot order is $c_ib_id_i$ and
$b\nsim c$ implies $a\nsim d$ as the knot order is $a_ic_ib_id_i$.
\endproof

In order to show the limit $abcd$ is a quadrisecant, we only have to
worry about whether the middle two intersection points $b$ and $c$
merge together. Pannwitz~\cite{Pann} proved that generic polygonal knots have a
quadrisecant and Kuperberg~\cite{Kup} proved that all (nontrivial
tame) knots have quadrisecants. He did this by using the limit argument
outlined above. In order 
to show the limit was a quadrisecant, he introduced the notion of
essential secants and quadrisecants (which he called ``topologically
nontrivial''). As we have encountered the same problem, we also use
the notion of essential, but redefine and extend it to suit our
purposes. We start by defining when an arc of a knot is essential,
capturing part of the knottedness of~$K$. Generically, a knot $K$
together with a secant segment $S=\overline{ab}$ forms a knotted
$\Theta$-graph in space. To adapt Kuperberg's definition, we consider
such knotted $\Theta$-graphs. A similar discussion is found in~\cite{Den,dds}.

\begin{figure}
\centerline{\begin{overpic}[scale=.5]{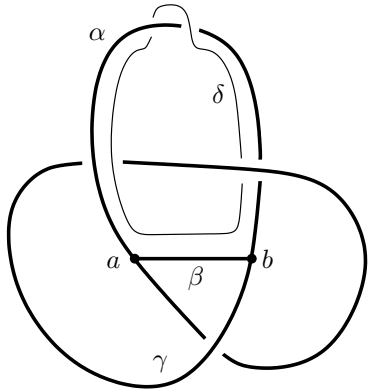}
\put(21.5,90){$\alpha$}
\put(26,31.5){$a$}
\put(66,31.5){$b$}
\put(47,27){$\beta$}
\put(38,6){$\gamma$}
\put(53.2,74){$\delta$}
\end{overpic}}
\caption[Essential arc of a knotted $\Theta$-graph.]{In the knotted
  $\Theta$-graph $\alpha\cup\beta\cup\gamma$, the ordered pair
  $(\alpha,\beta)$ is essential. To see this, we find the parallel
  $\delta=h(\alpha,\beta)$ to $\alpha\cup\beta$ which has linking
  number zero with $\alpha\cup\gamma$ and note that is is
  homotopically nontrivial in the knot complement $\R^3\setm(\alpha
  \cup\gamma)$. In this illustration, $\beta$ is the straight segment
  $\overline{ab}$, so we may equally say that the arc $\alpha$ of the
  knot $\alpha\cup\gamma$ is essential.}
\label{fig:essdef}
\end{figure}

\begin{definition}\label{def:esstheta}
  Suppose $\alpha$, $\beta$ and $\gamma$ are three disjoint simple
  arcs from~$a$ to~$b$, forming a knotted $\Theta$-graph.
  Then we say that the ordered pair $(\alpha,\beta)$ is \emph{inessential}
  if there is a disk~$D$ bounded by the knot $\alpha\cup\beta$ having no
  interior intersections with the knot $\alpha\cup\gamma$.
  (We allow self-intersections of~$D$, and interior intersections
  with $\beta$, as will be necessary if $\alpha\cup\beta$ is knotted.)
  
An equivalent definition is illustrated in \figr{essdef}:
  Let $X:=\R^3\setm (\alpha\cup\gamma)$, and let $\delta$ be a
  parallel curve to $\alpha\cup\beta$ in $X$. Here by \emph{parallel} we mean
  that $\alpha\cup\beta$ and $\delta$ cobound an annulus embedded in
  $X$. We choose $\delta$ so that it is homologically trivial in $X$
  (that is, so that $\delta$ has
  linking number zero with $\alpha\cup \gamma$). Let
  $h(\alpha,\beta)\in\pi_1(X)$ denote the (free) homotopy class of
  $\delta$. Then $(\alpha,\beta)$ is \emph{inessential} if
   $h(\alpha,\beta)$ is trivial.
We say that $(\alpha,\beta)$ is \emph{essential} if it is not
inessential. 
\end{definition}

This notion is clearly a topological invariant of the (ambient isotopy)
class of the knotted $\Theta$-graph. 
We apply this definition to arcs of a knot $K$
below. But first, some useful notation. Let $a,b \in K$. The arc from
$a$ to~$b$ following the orientation of 
the knot is denoted $\arc{ab}$. The arc from $b$ to $a$, $\arc{ba}$ is
similarly defined. The secant segment from $a$ to $b$ is denoted
$\overline{ab}$. We now define an essential arc of a knot following~\cite{dds},
rather than the stronger notion found in~\cite{Den}.

\begin{definition} \label{def:essarc}
  If $K$ is a knotted curve and $a,b\in K$, let $S=\segment{ab}$.
  We say $\arc{ab}$ is \emph{essential} if for every $\epsilon >0$ there exists some 
  $\epsilon$-perturbation of $S$ (with endpoints fixed) to a tame curve $S'$ such that   
  $K\cup S'$ forms an embedded $\Theta$ in which $(\arc{ab},S')$ is
  essential.
\end{definition}

Note that this definition is quite flexible. We will see that it ensures that the set of essential secants is closed in $S$. It also allows for the situation where $K$ intersects~$S$. We could allow perturbations only in this case, the combing arguments of~\cite{DS} show the resulting definitions are equivalent. Note also that we only require~$S'$ to be close to $S$ in the $C^0$ sense. Thus $S'$ could be locally knotted, however we only care about a homotopy class, not the isotopy class of~$h$.

In~\cite{ckks} it was shown that if $K$ is an unknot, then any arc
$\arc{ab}$ is inessential.  This follows immediately, because the
homology and homotopy groups of $X:=\R^3\setm K$ are equal for an
unknot, so any curve $\delta$ having linking number zero
with $K$ is homotopically trivial in $X$. Dehn's lemma \cite{Dehn,Pap,Rolf} is used to prove
a converse statement.

\begin{theorem} [\cite{dds}]
If $a,b\in K$ and both $\arc{ab}$ and $\arc{ba}$ are inessential,
then~$K$ is unknotted. \qed
\end{theorem}

\begin{figure}
\centerline{\begin{overpic}[scale=.4]{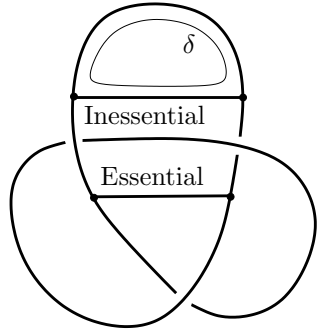}
\put(23,62.5){Inessential}
\put(28,44){Essential}
\put(53,84){$\delta$}
\end{overpic}}
\caption[Essential and inessential secants]{Examples of essential and
  inessential secants for the trefoil knot. The parallel curve
  $\delta$ is shown for the inessential secant.}
\label{fig:esecant}
\end{figure}

\begin{definition} 
  A secant $ab$ of $K$ is \emph{essential} if
  both subarcs $\arc{ab}$ and $\arc{ba}$ are essential. Otherwise it
  is \emph{inessential}. Let $ES$ be the \emph{set of essential secants
  in~$S$}.  \figr{esecant} illustrates an inessential and an essential
  secant of a trefoil knot.
\end{definition}

To call a quadrisecant $abcd$ essential, we could follow 
Kuperberg~\cite{Kup} and require that secants $ab$, $bc$ and $cd$
all be essential. However the notion of essential is needed only for
those secants whose endpoints are consecutive along the knot. The following definition is consistent with the more general definition found
in~\cite{dds}. 

\begin{figure}
\centerline{\begin{overpic}[scale=.4]{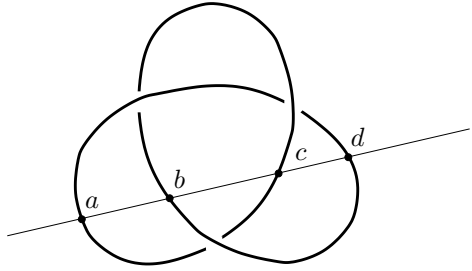}
\put(17,12.7){$a$}
\put(36,17){$b$}
\put(62,23){$c$}
\put(74,26){$d$}
\end{overpic}}
\caption[Essential alternating quadrisecant.]{Essential alternating
  quadrisecant $abcd$.}
\label{fig:equad}
\end{figure}

\begin{definition} \label{def:essnsec} An \emph{essential trisecant}
  $abc$ is a trisecant 
  which is essential in the second segment $bc$. An \emph{essential
  alternating quadrisecant} $abcd$ is an alternating quadrisecant
  which is essential in the second segment $bc$. (See \figr{equad}.)
  Trisecants and 
  alternating quadrisecants that are not essential are called \emph{inessential}. 
\end{definition}

We show in the proof of Main Theorem that when taking the limit of
essential secants, the points of the limit secant do not
combine to form one point, nor lie on a common straight subarc of $K$. Thus when
taking the limit of essential alternating quadrisecants, we may be
assured of an alternating quadrisecant. We want to prove that any
nontrivial generic polygonal knot has an \emph{essential} alternating
quadrisecant. Thus we need to understand the relationship between
essential trisecants and essential quadrisecants. 

\begin{definition}  Let $\mathcal{ET}$ be the set of essential
  trisecants in $K^3$, really, $\mathcal{ET}=\pi^{-1}_{23}(ES)\cap
  \mathcal{T}$.  Define $\mathcal{ET}^s$ and $\mathcal{ET}^d$ to be
  the sets of essential trisecants of same and different orderings in $K^3$.
Let $ET=\pi_{12}(\mathcal{ET})$ be the projection of the set of
  essential trisecants to the set of secants $S$ and similarly define
  $ET^s:=\pi_{12}(\mathcal{ET}^s)$ and $ET^d:=\pi_{12}(\mathcal{ET}^d)$.  
\end{definition}

\begin{lemma}\label{lem:essquadkey} Let $ab\in ET^s\cap ET^d$ in $S$. This
  means that there exists $c$ and~$d$ such that $abc\in \mathcal{ET}^d$ and
  $abd\in \mathcal{ET}^s$. Then either $abcd$ or $abdc$ is an essential
  alternating quadrisecant. 
\end{lemma}

\proof That $abcd$ or $abdc$ is an alternating quadrisecant follows
from \lem{alt}. To show it is essential, the midsegment $bc$,
respectively $bd$, must be essential.  The order along the quadrisecant
line is either $abcd$ or $abdc$. If it is $abcd$ then $bc$ is
essential as $abc\in \mathcal{ET}^d$. If it is $abdc$, then $bd$ is
essential as $abd\in \mathcal{ET}^s$.  
\endproof

Thus in order to show that a generic polygonal knot has at least one
essential alternating quadrisecant, it is sufficient to prove that
$ET^s$ and $ET^d$ have common points in $S$.

In the following Lemmas we describe the relationship between the set
of essential trisecants and the set of trisecants. We use this
relationship to completely understand the structure of $ET$ in $S$.

\begin{lemma}\label{lem:changeover} 
Let $K$ be a generic polygonal knotted curve. In $S\setm \pi_{13}(\mathcal{T})$, each
connected component consists of either inessential or essential
secants. Moreover, the set of essential secants is closed in~$S$,.
\end{lemma}

\proof  
Recall that $\pi_{13}:K^3\rightarrow K^2$ is projection onto the first
and third coordinates. \prop{allproj} showed that $\pi_{13}$ is
a piecewise smooth immersion of $\mathcal{T}$ in $S$ and
$\pi_{13}(\mathcal{T})$ intersects itself in double points.

\begin{figure}
\centerline{\begin{overpic}[scale=.4]{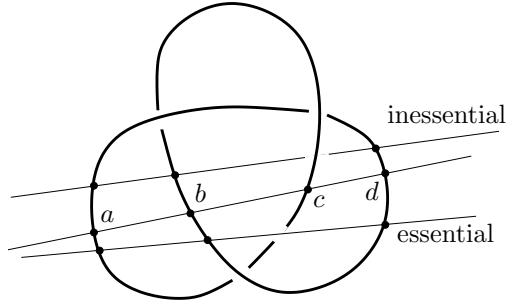}
\put(19,16){$a$}
\put(38,20){$b$}
\put(62,19){$c$}
\put(72.5,20.5){$d$}
\put(77,36){inessential}
\put(79,12){essential}
\end{overpic}}
\caption[Secant $bd$ changes from inessential to essential through
trisecant $bcd$]{Secant $bd$ changes from inessential to essential through
trisecant $bcd$ and so trisecant $abd$ changes from inessential to
essential through quadrisecant $abcd$.}\label{fig:etri}
\end{figure}

We wish to show that in $S$, secants change from inessential to essential
\emph{only} when there is a trisecant $t\in\pi_{13}(\mathcal{T})$. Let
$ab$ be an inessential (essential) secant. By perturbing $a,b\in K$ we
can get to all nearby secants in $S$. As long as the chord
$\overline{ab}$ never touches $K$ during the perturbation then the
topological type of the knotted $\Theta$-graph has not changed.
Thus all secants near secant $ab$ will have the same topological
type in the knotted $\Theta$-graph \emph{unless} $K$ intersects
$\overline{ab}$. In this case, $ab\in\pi_{13}(\mathcal{T})$. Thus
secants change from inessential to essential \emph{only} when there is
a trisecant $t\in \pi_{13}(\mathcal{T})$ and each connected
component of $S\setm \pi_{13}(\mathcal{T})$ consists of either
essential or inessential secants. (For example, \figr{etri} shows a
trisecant $abd$ changing from inessential to essential through
quadrisecant $abcd$. Here secant~$bd$ changes from
inessential to essential via trisecant $bcd$.)

Moreover, $ES$ is closed in $S$. Take a sequence of essential secants
$\{s_i\}_{i=1}^{\infty}$, and suppose $s_i\rightarrow s$ where
$s\in\pi_{13}(\mathcal{T})$. First assume $s_i\notin
\pi_{13}(\mathcal{T})$. Then for any $\epsilon >0$, far enough along
the sequence each~$s_i$ 
is an $\epsilon$-perturbation of~$s$. (By being within~$\epsilon/2$ of an~$\epsilon/2$-perturbation of~$s$.) Thus by definition~$s$ is  
essential. Now assume $s_i\in \pi_{13}(\mathcal{T})$. Then
for each~$s_i$ there is a sequence $s_i^j\rightarrow s_i$, where~$s_i^j$ are essential secants and $s_i^j\notin
\pi_{13}(\mathcal{T})$. For each~$i$, there is a~$j_i$ such that the sequence~$\{s_i^{j_i}\}_i$ of essential secants (not in
$\pi_{13}(\mathcal{T})$) converges to~$s$. Again, for any
$\epsilon >0$, far enough along the sequence each $s_i^{j_i}$ is an
$\epsilon$-perturbation of~$s$ and thus by definition~$s$ is essential. 
\endproof

\begin{corollary}
In $S$, intervals of $\pi_{13}(\mathcal{T})$ between two connected
components of essential secants are essential. Intervals of
$\pi_{13}(\mathcal{T})$ between a component of essential and a
component of inessential secants are essential.
\end{corollary}

\proof Apply \lem{changeover}. \endproof

Note that \lem{changeover} says nothing about intervals of
$\pi_{13}(\mathcal{T})$ between two connected components of
inessential secants. For non-generic knots it is possible to have such
isolated essential secants, but this case does not matter for us.

We use this information to describe exactly how the sets of trisecants
and essential trisecants are related. Let $\overline{\mathcal{ET}}$
denote the closure of $\mathcal{ET}$ in $K^3$ and let
$\overline{ET}=\pi_{12}(\overline{\mathcal{ET}})$. In a similar way
$\overline{ET}^s=\pi_{12}(\overline{\mathcal{ET}}^s)$ and
$\overline{ET}^d=\pi_{12}(\overline{\mathcal{ET}}^d)$.

\begin{lemma}\label{lem:essclosed} Let $K$ be a generic polygonal
  knot. Then $\mathcal{ET}$ is a closed subset of~$\mathcal{T}$ and
  $\overline{\mathcal{ET}}$ is a finite union of closed circles and
  closed intervals in $\overline{\mathcal{T}}$. 
\end{lemma}

\proof
Trisecant $abc$ is inessential when the second segment $bc$ is
essential. By definition, $\mathcal{ET}=\pi^{-1}_{23}(ES)\cap
\mathcal{T}$. Now $ES$ is closed in $S$ and $\pi^{-1}_{23}(ES)$ is
closed in~$K^3$. Thus $\mathcal{ET}$ is a closed subset of $\mathcal{T}$.
By \lem{changeover}, connected components of $S\setm
\pi_{13}(\mathcal{T})$ are either inessential or essential
secants. By examining where $\pi_{23}(\mathcal{T})$ lies in relation
to the connected components of essential secants, we may
determine which trisecants are  essential or inessential.  

Let $e_i\times e_j$ be the subset of $S$ consisting of all secants with points lying on edges
$e_i$ and $e_j$ in that order. In  $e_i\times e_j$, 
$\pi_{13}(\mathcal{T})$ and $\pi_{23}(\mathcal{T})$ are intervals of trisecants. By modifying the arguments found in \lem{onequad}, these intervals either do not
intersect or intersect in one or two (distinct) points. It is only at these points of intersection, that $\pi_{23}(\mathcal{T})$ can change from inessential to essential.
Thus in each $e_i\times e_j \subset K^2$, $\pi_{23}(\mathcal{T})$ either
remains inessential, remains essential or changes
from inessential to essential at most \emph{twice}.
Hence $\overline{\mathcal{ET}}$ is a \emph{finite} union of closed circles and
closed intervals in $\overline{\mathcal{T}}$. \endproof

\begin{lemma}\label{lem:immersedess} Let $K$ be a generic polygonal knotted
  curve. Projection $\pi_{12}$ is a piecewise smooth immersion of
  $\mathcal{ET}$ to $S$. In $\overline{T}$, $\overline{ET}$ is a finite
union of closed intervals and closed circles and in $S$, $ET$
  intersects itself in double points. 
\end{lemma}

\proof As ${\mathcal{ET}}$ is closed in ${\mathcal{T}}$, \lem{immersed} shows projection $\pi_{12}$ is a piecewise smooth immersion of  $\mathcal{ET}$ to $S$. Thus $ET$ inherits the properties of $T$
found in Lemmas~\ref{lem:onequad} and \ref{lem:doublepoints}. As 
$\overline{\mathcal{ET}}$ is also closed in $\overline{\mathcal{T}}$, $\pi_{12}(\overline{\mathcal{ET}})$ is a finite union of closed
intervals and circles.
 \endproof

\begin{lemma}\label{lem:essdist} 
 Let $K$ be a generic polygonal knotted curve. In $S$,
  $d_S(\overline{ET}^d,\tilde{\Delta}_-)\geq h$ and
  $d_S(\overline{ET}^s,\tilde{\Delta}_+)\geq 
  h$, where $h$ is the minimum edge length of $K$.
\end{lemma}

\proof As $\overline{ET}$ is a subset of $\overline{T}$, \lem{dist}
gives the result.  
\endproof

Thus $\overline{ET}^s$ and $\overline{ET}^d$ can only intersect in $S=K^2\setm\tilde{\Delta}$. As
$ET$ is closed in~$T$, the boundary of $ET$ will contain points other
than $\overline{ET}\setm ET$ and can be found anywhere in
$\overline{S}= S\cup\tilde{\Delta}_+ \cup \tilde{\Delta}_- $. The content of the next 
lemma is that the points in $\overline{ET}\setm ET$ must occur on
$\Delta_+$ or $\Delta_-$.

\begin{lemma}\label{lem:bdryet}
In $\overline{S}$, $\pi_{12}(\overline{\mathcal{ET}} \cap \bdry \mathcal{T})
\subset \Delta_-$ or $\Delta_+$. 
\end{lemma}

\begin{figure}
\centerline{\begin{overpic}[scale=.4]{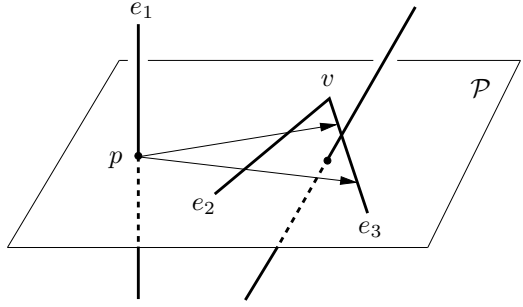}
\put(20,27){$p$}
\put(24,56){$e_1$}
\put(36,18){$e_2$}
\put(68.3,13.5){$e_3$}
\put(61,42){$v$}
\put(90,40){$\mathcal{P}$}
\end{overpic}}
\caption[Trisecants whose last two points converge cannot be essential
trisecants.]{The interval of trisecants through $pe_2e_3$ must
  eventually be inessential as it ends at degenerate trisecant
  $pvv$. Close to vertex $v$ the second segment through $e_2e_3$ is
  inessential. Trisecants through $e_3e_2p$ may be essential until 
  they end at degenerate trisecant $vvp$ as the second segment $e_2p$
  may be essential.}
\label{fig:inesstri}
\end{figure}

\proof
In $\overline{S}$, we know that $\bdry T=\pi_{12}(\bdry\mathcal{T})$ occurs 
when an interval of adjacent trisecants ends on a degenerate trisecant $pvv$ or $vvp$ (see
\lem{allint}). From discussion in \sect{tri}, $\bdry T$ occurs in $S$ when the degenerate trisecant is of the form~$pvv$. Let the points of the trisecant lie on edges $e_1$,~$e_2$ and~$e_3$ in that order and let ~$e_2$ and~$e_3$ be adjacent edges with common vertex~$v$. Let~$e_1$ intersect the plane~$\mathcal{P}$
spanned by~$e_2$ and~$e_3$ at point~$p$. See
\figr{inesstri}. Trisecants with points on~$pe_2e_3$ must eventually be inessential. This is because the second segment
through~$e_2$ and~$e_3$ must eventually be inessential. As the trisecant
converges to~$pvv$, there will be an embedded disk in
$\mathcal{P}$ spanned by parts of $e_2$, $e_3$ and the part of the
trisecant line between $e_2$ and $e_3$.

On the other hand, the part of $\bdry T$ that occurs on $\Delta_+$ or
$\Delta_-$ might contain $\overline{ET}$. Use the same set up as in
the previous paragraph and in \figr{inesstri}, but assume the points
of the trisecant lie on $e_3e_2e_1$ in that order. Here the interval of adjacent trisecants end on degenerate trisecant $vvp$. The second segment $e_2p$
of the trisecant may be essential. Thus there is nothing to prevent
an interval of essential trisecants ending at degenerate trisecant
$vvp$ and $\pi_{12}(vvp)=vv\in \Delta_+\cup\Delta_-$.
\endproof

Thus unlike $T$, $ET$ (and also $ET^s$, $ET^d$) are closed in $S$. Moreover, as points of $\overline{ET}\setm ET$ must occur on $\Delta_-$ or $\Delta_+$, we need not be concerned about degenerate trisecants when considering $ET^s\cap ET^d$. Hence there is no need for a result similar to \lem{intint} in order to show that generic polygonal knotted curves have an essential alternating quadrisecant. We have all the information needed. 



\section{Essential Alternating Quadrisecants for Generic Polygonal
Knots}\label{sect:quadgenpoly}  

\lem{essquadkey} and \lem{bdryet} showed that $ET^s \cap {ET}^d \neq\emptyset$ in ${S}$ is all that is needed to prove that an essential alternating quadrisecant exists for generic polygonal
knotted curves. 

The underlying motivation for the proofs in this section is that idea
that every curve going from left to right across a square is
intersected by every curve going from the top to the bottom. In our
case, we don't have a square, but an annulus. This may be thought of
as a square with the left and right edges identified. Thus every curve
winding once around the annulus is 
intersected by every curve going across (from top to bottom).

In \lem{esspann} we show that all curves winding once around $S$ (defined
below) must be intersected by $ET$. Thus intuitively $ET$ ``goes across'' the
annulus $S$. This together with the fact that ${ET}^s$ avoids the top of
$S$ and ${ET}^d$ avoids the bottom of~$S$, implies that ${ET}^s\cap
{ET}^d\neq\emptyset$ in $S$.

\begin{figure}
\centerline{\begin{overpic}[scale=.4]{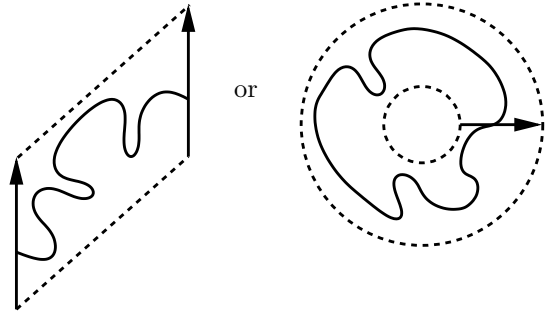}
\put(42,40){or}
\end{overpic}}
\caption[Closed curve winding once around $S$.]{A curve winding
  once around $S$ is a simple closed curve homotopic to the generator $1$ in
  $\pi_1(S)\equiv \mathbb{Z}$.}\label{fig:curve}
\end{figure}

\begin{definition} A closed curve winding once \emph{around} $S$ is a
  simple (continuous) curve $\alpha:[0,1]\rightarrow S$ which is homotopic to 1 in
  $\pi_1(S)\equiv \mathbb{Z}$. See \figr{curve}.
\end{definition}

\begin{lemma}\label{lem:esspann} (Pannwitz) Let $K$ be a nontrivial generic
  polygonal knotted curve. Every closed curve winding once around
  $S$ has nonempty intersection with the set of essential trisecants. 
\end{lemma}

\proof 
This proof contains ideas from each of \cite{Pann, Schm, Kup}.
Assume, by way of contradiction, that there is a curve~$\alpha$
winding once around $S$ with \emph{empty} intersection with the set of
essential trisecants. As $\alpha$ avoids the set of essential trisecants,
each point of $\alpha$ is the first two points of a secant or an
\emph{inessential} trisecant. To simplify the argument we also assume
that the curve $\alpha$ intersects $\pi_{12}(\mathcal{T}\setm\m{ET})$ transversely
away from self-intersections. As the set of essential secants is
closed in ${S}$, curves avoiding~$ET$ are in an open
set. Any curve which is not transverse to $\pi_{12}(\mathcal{T}\setm\m{ET})$ has another curve arbitrarily close to it which is. Also any curve which passes through a self-intersection of  $\pi_{12}(\mathcal{T}\setm\m{ET})$
has another curve arbitrarily close to it which avoids that self-intersection. Once we have
proved the theorem for curves which intersect
$\pi_{12}(\mathcal{T}\setm\m{ET})$ transversely away from self-intersections, we
have proved it for \emph{all} curves winding once around~$S$.

\begin{figure}
\centerline{\begin{overpic}[scale=.4]{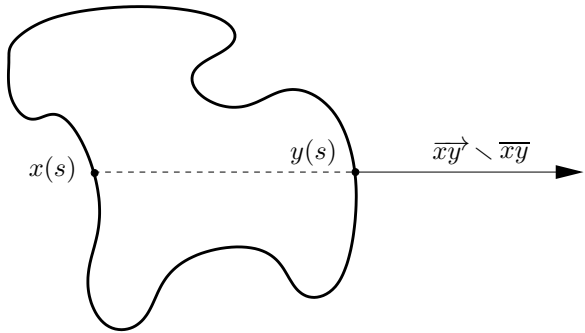}
\put(4,27){$x(s)$}
\put(49,30){$y(s)$}
\put(73,30){$\overrightarrow{xy}\setm \overline{xy}$}
\end{overpic}}
\caption[Rays from a knot.]{An example of a ray $\overrightarrow{xy}\setm\overline{xy}$. The union of such rays and the
  point at infinity 
  form a disk with the knot as its boundary.}\label{fig:ray}
\end{figure}

Let us first assume each point of $\alpha$ is a secant. We will deal with
the case where a point of $\alpha$ is the first two points of an inessential trisecant
later.
Let~$\alpha=(x(s),y(s))$  where we have parameterized $[0,1]$ with respect
to arclength. Construct all (geodesic) rays
$\overrightarrow{xy}\setm \overline{xy}$. This is the
collection of rays that start at~$y(s)$ (the second
point of the secant) and extend to infinity in the direction
of the vector from $x(s)$ to $y(s)$. Such a ray is illustrated in \figr{ray}. 

With the point at infinity, the union of the rays forms a
spanning disk~$D$ with~$K$ as its boundary. This is because of the
continuity of the curve $\alpha$ and the fact that $\alpha$ is a closed curve
homotopic the generator $1$ in $\pi_1(S)$. Thus both $x(s)$ and $y(s)$ have degree $1$ on the annulus $S$, and hence traverse the knot once. The assumption
that each point on $\alpha$ is a secant means that the rays do
not intersect the knot again. However, the disk may have self-intersections. In particular there may be self-intersections on the boundary of the
disk as $y(s)$ is not necessarily strictly monotonic around $K$. 

The aim is to apply Dehn's Lemma. To do
so, we must find a disk bounded by~
$K$ which is the image of a PL-map and whose
self-intersections do not occur on the boundary. We do this by
altering the disk $D$ constructed above. Construct an
$\epsilon$-neighborhood~$N$ around~$K$. (This is the set of
all points within $\epsilon$ of the knot.) We take $\epsilon$ small
enough so that this a regular neighborhood. If necessary perturb~$N$ 
slightly so that $D$ intersects the boundary of $N$, $\bdry N$, transversely.
The disk $D$ intersects~$\bdry{N}$ in
several kinds of closed curves. Most are homologous to~0 on $\bdry N$,
but one will be homologous to a curve $\gamma$ of homotopy class $(1,m)$. We wish to replace $D$ with a disk that only intersects~$\bdry N$ in $\gamma$. 

The curves homologous to $0$ on $\bdry {N}$ are disjoint, thus
abstractly, these 
intersection curves are many families of nested circles. To remove these curves, take the outermost intersection circle of a family of nested
circles which has some points in its interior lying inside~$N$. The
circle bounds a disk $A$ on $D$ and a disk $B$ on $\bdry N$. Replace $A$ by
$B$ and push~$B$ slightly off $\bdry N$ to simplify the
intersection. Repeating this process for all other families of nested
circles eliminates all curves homologous to $0$. Call the new disk so
constructed~$D'$ and note that $D'$ intersects $\bdry{N}$ only in the curve~$\gamma$.  

We wish to replace the part of $D'$ outside $\bdry{N}$ with a PL disk.
Approximate all of $D'$ by a PL disk and use this to approximate $\gamma$
very closely by a PL curve $\gamma_0$ outside 
$\bdry N$. (Note $\gamma_0$ is of the same homology class as $\gamma$
in the solid 
toroidal shell they both lie on.) Remove the part of the
PL disk inbetween $\gamma_0$ and $K$ (this is mostly inside $N$).  Thus we
have a PL-disk with $\gamma_0$ as its boundary. We would like to join
$K$ to $\gamma_0$ with an embedded PL collar, but can't as $\gamma_0$
might have self intersections. Now $\gamma_0$ is homologous to a
curve of homotopy class $(1,m)$.  Put a smaller regular $\epsilon$-neighborhood
$N_1$ inside $N$. On $\bdry N_1$ put an embedded $(1,m)$ curve denoted
$\gamma_1$. Thus $\gamma_0$ and $\gamma_1$ live on the boundary of a
solid toroidal shell where they are homologous. We use a PL-homotopy
in this shell to join $\gamma_0$ and~$\gamma_1$ and then join
$\gamma_1$ to~$K$ by an embedded PL-collar.  This creates a new PL-disk
with boundary $K$ and with no boundary intersections. By Dehn's
Lemma, we may replace this new disk by an embedded one, giving a
contradiction to knottedness.
 
\begin{figure}
\centerline{\includegraphics{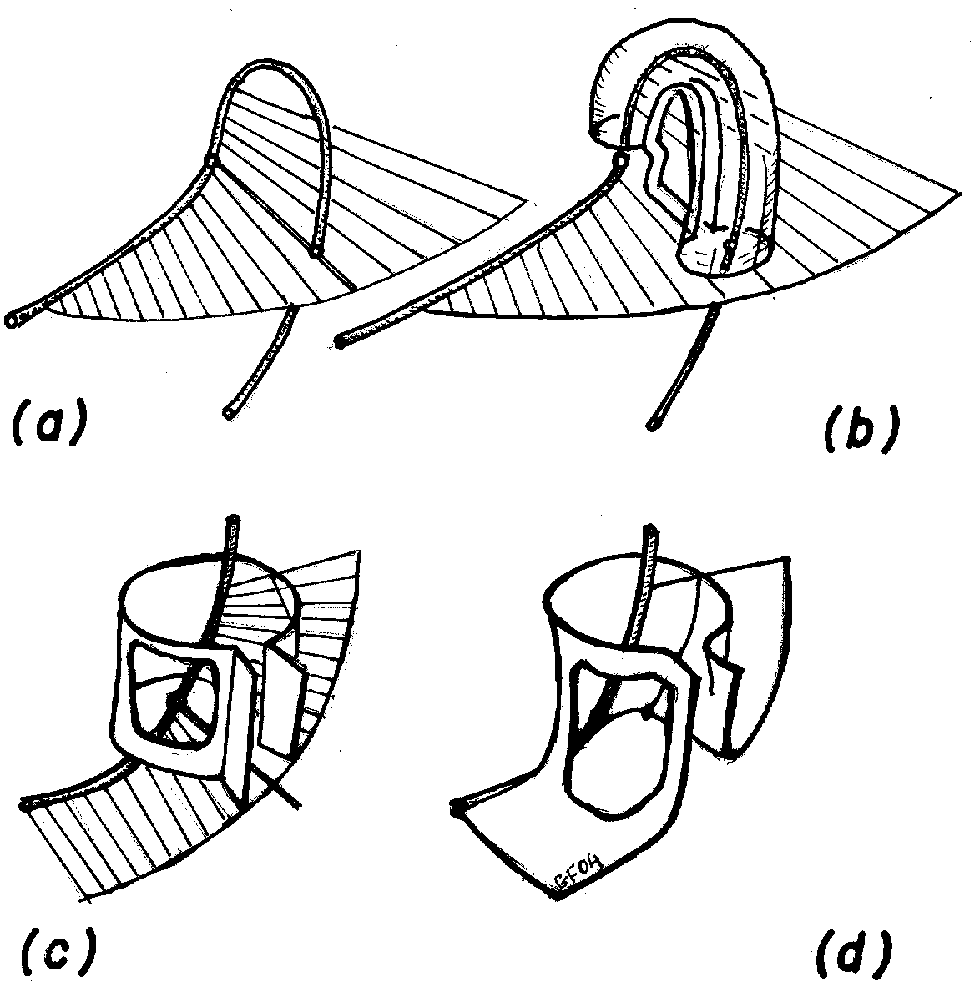}}
\caption{Figure illustrating the surgery described in \lem{esspann}.}
\label{fig:kshmoosh}
\end{figure}

Now on the original curve $\alpha$ there may be a point which is the first
two points of an inessential trisecant $abc$. We assumed that $\alpha$ intersects $\pi_{12}(\mathcal{T})$ transversely away from self-intersections, thus such a
point is isolated. In order to apply Dehn's lemma, we perform surgery on the disk $D$ so that the third point~$c$ of the trisecant no longer intersects the
disk. As the second segment $\overline{bc}$ of trisecant $abc$ is
inessential, segment $\overline{bc}$ and one arc of
the knot spans a disk $\mathcal{D}$ whose interior avoids the
knot. (See \figr{kshmoosh}(a)\footnote{Thanks to G. Francis for permission to use \figr{kshmoosh} and \figr{kupdisk}.} .) The spanning disk $\mathcal{D}$
intersects the disk $D$ in the segment $\overline{bc}$. Take a small
embedded $\epsilon$-neighborhood of the knot, this intersects both the
original disk $D$ and the spanning disk $\mathcal{D}$. Now make two
copies of the spanning disk $\mathcal{D}$ and move them apart. Remove
the parts of the copies of $\mathcal{D}$ inside the
$\epsilon$-neighborhood. Also remove the parts of $D$ inside the two
copies of $\mathcal{D}$ (this is a small strip about the inessential
segment $\overline{bc}$). (See \figr{kshmoosh}(b).) After smoothing
things out, the original disk $D$ has been altered so that it does not
intersect the third point of the trisecant. However there is a self-intersection 
near the second point $b$ of the trisecant. This is
illustrated in \figr{kshmoosh}(c) and~(d). This surgery may be
repeated for any other points of the curve~$\alpha$ which are the first two points of
an inessential trisecant. Thus we obtain a disk $D$ whose interior
avoids $K$ and we apply Dehn's Lemma as before to get a contradiction
to knottedness.  
\endproof

\begin{proposition}\label{prop:homology}
Let $A$ and $B$ be two closed subsets of the annulus
${S}$, such that $A$ lies outside some neighborhood of
$\Delta_+$ and $B$ lies outside some neighborhood of $\Delta_-$. If $A\cap
B=\emptyset$, then there is a curve winding once around $S$ avoiding
$A\cup B$.
\end{proposition}

\proof
Look at the Mayer-Vietoris sequence:
$$ \dots\rightarrow H_1(S\setm(A\cup B))\rightarrow H_1(S\setm
A)\oplus H_1(S\setm B) \rightarrow H_1(S)\rightarrow\dots $$

Using the assumption that $A$ lies outside some neighborhood $U$ of
$\Delta_+$,  construct a path~$\alpha$ in $U\subset S\setm A$ which
winds once around~$S$. Similarly, as $B$ lies outside some
neighborhood $V$ of $\Delta_-$, construct a path~$\beta$ in
$V\subset S\setm B$ which winds once around~$S$ with reverse
orientation. These paths represent 
homology classes in $H_1(S\setm A)$ and $H_1(S\setm B)$ and the image
of these classes in $H_1(S)$ is homologous to~+1 and~-1
respectively. The map $f: H_1(S\setm 
A)\oplus H_1(S\setm B) \rightarrow H_1(S)$ takes a pair
$(p,q)$ to the sum of the image of~$p$ and the image of~$q$ in
$H_1(S)$. Thus the image 
of $(\alpha,\beta)$ under~$f$ is $0$. Now the map $g: H_1(S\setm(A\cup
B))\rightarrow H_1(S\setm A)\oplus H_1(S\setm B)$ takes a class
$\gamma$ to the pair $(p,-q)$ where $p$ is the image of
$\gamma$ in $H_1(S\setm A)$ and $q$ is the image of $\gamma$ in
$H_1(S\setm B)$.  Therefore, by exactness of the sequence at
$H_1(S\setm A)\oplus 
H_1(S\setm B)$, there is a $\gamma$ in $H_1(S\setm(A\cup B))$ which maps to
$(\alpha,\beta)$. This $\gamma$ may be represented by a path in
$S\setm(A\cup B)$. As $\gamma$ maps to $\alpha$, $\gamma$ winds once
around $S$ and avoids~$A\cup B$.\endproof

\begin{theorem}\label{thm:equad}  Every nontrivial generic polygonal knotted
  curve in $\R^3$ has an essential alternating quadrisecant.
\end{theorem}

\proof  \lem{essdist}, together with $ET$ closed in $S$, shows \prop{homology} may be applied to ${ET}^s$ and ${ET}^d$ in ${S}$. If ${ET}^s\cap {ET}^d =\emptyset$ in ${S}$, then there exists a path winding once
around $S$ avoiding $ET={ET}^s\cup {ET}^d$. This is a
contradiction to \lem{esspann} (Pannwitz). Hence in $S$,
${ET}^s\cap{ET}^d\neq\emptyset$, 
which by \lem{essquadkey} shows that there is at least one essential
alternating quadrisecant.
\endproof

\section{Main Result}\label{sect:mainresult}

We now extend \thm{equad} by using a limit
argument to show that essential alternating quadrisecants exist for \emph{all}
nontrivial tame knots $K$ in $\R^3$. First, we define the kind of
closed curves and the kind of convergence that we need for the limit
to make sense. 

Milnor~\cite{Mil} defined the total curvature of an arbitrary curve
as the supremal total curvature of inscribed polygons.  Curves of finite
total curvature are exactly those rectifiable curves whose unit
tangent vectors have bounded variation. Any curve of finite total
curvature has well-defined one-sided 
tangent vectors everywhere; these are equal and opposite except at
countably many \emph{corners}. (See also discussion of curves of finite total curvature in CITE.) When considering limits of knots of finite total curvature, we will need to have control over the behavior of tangent vectors. Because knots of finite total curvature have corners, we now define a notion of convergence similar to, but more general than~$C^1$-convergence.

\begin{definition} 
Suppose $\{K_i\}_{i=1}^\infty$ is a sequence of (closed) curves of finite total curvature. 
Then we say $K_i$ \emph{converges in a $C^1$ sense} to a limit
curve $K$ if there is a set of parameterizations of $K_i$ such that for
all $\epsilon >0$ there is an $I$ such that once $i>I$ then for all $t$,
$K_i(t)$ is within $\epsilon$ of~$K(t)$. Also for all $t$, left
and right unit tangent vectors of $K_i(t)$ and $K(t)$ are within
$\epsilon$ of each other. 
\end{definition}

\begin{theorem}[Main Theorem]\label{thm:mainthma}
Every knotted curve in $\R^3$ has an alternating quadrisecant.
\end{theorem}

\proof This uses ideas from \cite{Kup,Schm} that have been altered and extended to
suit this case.
Given a knotted curve~$K$, since it is tame, there is a homeomorphism
$h:\mathbb{R}^3 \rightarrow \mathbb{R}^3$ such that~$h(K)$ is
smooth. Furthermore~$h$ may be chosen so that~$h$ is smooth on
$\mathbb{R}^3 \setminus K$ (see \cite{Kup}). 

Use \lem{gpconverge} to construct a sequence of generic polygonal
knots $\{K_i\}_{i=1}^{\infty}$ such that $K_i$ converges to $K$ and
with the following properties:
\item{(1)} $h(K_i)\cap h(K)=\emptyset$
\item{(2)} $h(K_i)$ converges in a $C^1$ sense to $h(K)$.

By construction each generic polygonal~$K_i$ has the same isotopy type
as~$K$, thus each~$K_i$ is nontrivial. By \thm{equad} we know that
each~$K_i$ has an essential alternating
quadrisecant~$a_ib_ic_id_i$. By picking a subsequence if necessary we
may assume that the original sequence converges to~$abcd\in K^4$. We
also assume the order of~$a_ib_ic_id_i$ along the knot is~$a_ic_ib_id_i$. From \lem{converge}, it is sufficient to show that~$b$
and~$c$ do not lie on the same straight subarc of~$K$ in order to
show~$abcd$ is a quadrisecant. Thus we examine the limiting behavior of $b_ic_i$. Recall that each (generic) $K_i$ has no quintisecants,  hence $K_i\cap\overline{b_ic_i}=\{b_i,c_i\}$~only.

\begin{figure}
\centerline{\includegraphics{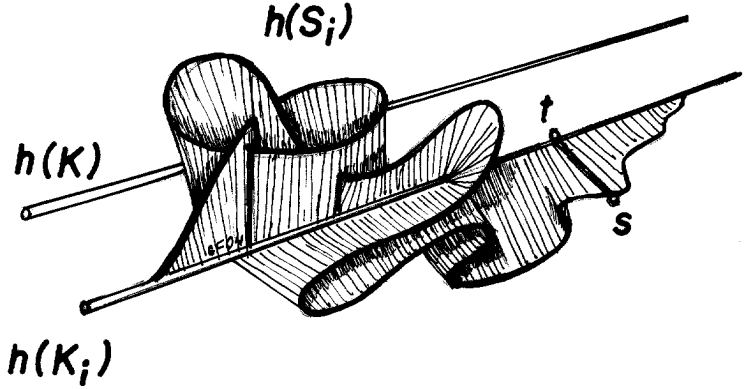}}
\caption[Disk constructed in the proof of \thm{mainthma}.]{Disk
  constructed in the proof of \thm{mainthma}. The boundary of $D$ is
  $h(S_i)\cup h(\arc{S_i})$ and $D$ is constructed as the union of
  line segments from $h(S_i)$ to $h(\arc{S_i})$.}
\label{fig:kupdisk}
\end{figure}

Let $S_i$ be the essential segment $\overline{b_ic_i}$, let $\arc{S_i}=\arc{b_i c_i}$ and suppose, by way of
contradiction, that~$b$ and~$c$ lie on the same straight subarc of~$K$. Let $S=\overline{bc}$. Then $\arc{bc}=S$ (thus is contained in $K$) and 
$h(S)=h(\arc{bc})\subset h(K)$. Take a small neighborhood of
$h(\arc{bc})$ that lies inside the embedded normal tubular neighborhood
about $h(K)$. As $h(\arc{bc})=h(S)$, far enough down the sequence
$h(S_i)$ and $h(\arc{S_i})$ both lie inside this
neighborhood.
For each point $s\in h(S_i)$
take the unique point $t\in h(K_i)$ lying in the same normal disk to~$h(K)$. Take the union of line segments from~$s$ to~$t$, as
illustrated in \figr{kupdisk}.
This forms a spanning disk~$D$ of $h(\arc{S_i})\cup h(S_i)$. But $h(K_i)$ does not intersect  $D^\circ$, the interior of $D$, else $h(\arc{S_i})$ is not a section. (In \figr{kupdisk} we see~$h(S_i)$
intersects~$D^\circ$ and $D^\circ$ has intersections, as allowed in the definition of inessential.)
 Thus for the knotted theta
$\Theta_i=h(K_i)\cup h(S_i)$, the pair $(h(\arc{S_i}),h(S_i))$ is
inessential. As essential is a topological notion for a knotted theta, the pair $(\arc{S_i},S_i)$ is
inessential and hence secant~$b_ic_i$ is inessential, a contradiction. \endproof

All that remains is to construct the desired sequence of generic
polygonal knots. 

\begin{lemma}\label{lem:gpconverge}
Let $K$ be a knotted curve and let $h$ be the homeomorphism
$h:\mathbb{R}^3 \rightarrow  
\mathbb{R}^3$ described in \cite{Kup} such that $h(K)$ is smooth and
$h$ is smooth on $\mathbb{R}^3 \setminus K$. Then there is a sequence of generic polygonal
knots $\{K_i\}_{i=1}^{\infty}$ such that~$K_i$ converges to~$K$
with the following properties:
\item{(1)} $h(K_i)\cap h(K)=\emptyset$
\item{(2)} $h(K_i)$ converges in a $C^1$ sense to $h(K)$.
\end{lemma}

\proof
As $h(K)$ is smooth, there is a sequence of smooth knots
$\{K^s_i\}_{i=1}^\infty$ such that $h(K^s_i)\cap h(K)=\emptyset$
and $h(K^s_i)$ converges in $C^1$ to $h(K)$. 
For each smooth $K_i^s$, there is a family of
inscribed polygons which are ambient isotopic to $K_i^s$ (see for instance \cite{CF}). We now show that these
polygonal knots may be perturbed to give a sequence
$\{K_i\}_{i=1}^{\infty}$ of generic 
polygonal knots whose images under $h$ are ``close'' to each $h(K_i^s)$ and hence converge in a~$C^1$ sense to $h(K)$.

By assumption each $h(K_i^s)$ is disjoint from $h(K)$ and so lies outside an open normal embedded tube of radius $\epsilon_i$ about $h(K)$. In fact $h(K_i^s)$ can be chosen to lie in a toroidal shell between two tubes of radius $\epsilon_i$ and $2\epsilon_i$ respectively about $h(K)$.
The map $h$ is only smooth on the complement of~$K$ and its
derivatives may become unbounded at $K$ (if $K$ is not smooth for
example). However, as a continuous map is bounded on a compact set, the derivatives of~$h$ are bounded in the toroidal shell that $h(K_i^s)$ lies inside.
So the ratio of the maximum to minimum values of both the first and the second derivatives of $h$ is bounded by some constant. Because these derivatives are bounded, it is a straightforward matter to find a generic polygonal~$K_i$ sufficiently close to $K^s_i$, so that if two corresponding points and tangent vectors are very close to one another in distance and angle respectively, then the image of the points and vectors under $h$ will be very close as well. As $h(K^s_i)$ converges in $C^1$ to $h(K)$, then $h(K_i)$ converges in a $C^1$ sense to $h(K)$. Thus
$\{K_i\}_{i=1}^\infty$ is the desired sequence of generic polygonal knots.
 \endproof

In essence \thm{mainthma} showed that the limit of essential secants is still a secant. In \cite{dds} we showed (with some extra assumptions) that for a nontrivial $C^{1,1}$ knot a limit of essential secants is still essential. In fact this result can be generalized to nontrivial knots of finite total curvature in one of two ways. The first is to generalize the arguments found in \cite{dds}. The second is to recall that being essential is a topological property of a knotted $\Theta$-graph. So in order to show that a limit of essential secants remains essential, it is sufficient to show that nearby knotted $\Theta$s are isotopic. This approach has been pursued in another paper \cite{DS}. We also conjecture that any nontrival (tame) knot has essential quadrisecants, but in the meantime we claim the following.

\begin{corollary}\label{cor:altquadess}
Every nontrivial knot of finite total curvature in $\R^3$ has an essential alternating quadrisecant.
\end{corollary}
\begin{proof}
Given a knot of finite total curvature, there is a family of inscribed polygons which are ambient isotopic to it (see for instance \cite{AB}). By perturbing these slightly, we find a sequence of generic polygonal knots converging in a $C^1$ sense to $K$. By \thm{mainthma} these have essential alternating quadrisecants. Some subsequence of these quadrisecants converges to an alternating quadrisecant, which is essential (as proved in \cite{DS}).
\end{proof}

\section{Corollaries to the Main Theorem}\label{sect:corollaries}

The Main Theorem has two immediate applications. It is used to give
alternate proofs to two previously known theorems about the geometry
of knotted curves. The existence of alternating
quadrisecants for knotted curves captures the intuition that a space
curve must loop around at least twice to become knotted. This
intuition is also reflected in results about total curvature and
second hull of knotted curves.

Recall that for polygonal closed curves the total curvature is the sum of the exterior angles and Milnor \cite{Mil} defined the total curvature of an arbitrary curve as the supremal total curvature of inscribed polygons. In 1929,
M.~Fenchel~\cite{Fen} proved that the total curvature of a closed
curve in~$\R^3$ is greater than or equal to~$2\pi$, equality holding
only for plane convex curves. In~1947, K.~Borsuk~\cite{Bor} extended this result to~$\R^n$ and conjectured the following:

\begin{theorem} A knotted curve in $\R^3$ has total curvature
  greater than $4\pi$.
\end{theorem}

This result was first proved around 1949 by both Milnor~\cite{Mil} and
I.~F\'ary~\cite{Fary}. It has since become known as the F\'ary-Milnor
theorem. In his proof, Milnor
used the idea of bridge number\footnote{The 
  bridge number $b(K)$ of a knot $K$ is the minimum number of bridges
  (or overpasses) occurring in a diagram of the knot, where the
  minimum is taken over all possible diagrams of $K$.}, making the
observation that for a knotted curve, there are planes in every
direction which cut the knot at least four times. More recently
in \cite{ckks}, the theorem was proved using the existence of second
hull for a knotted curve (defined later). Here we give a new proof using
the existence of alternating quadrisecants.
 
\proof  A  knotted curve $K$ has an alternating quadrisecant by \thm{mainthma}. An alternating quadrisecant is an inscribed
polygon in $K$ with total curvature~$4\pi$. By definition, its total curvature is less than or equal to the total curvature of the knotted curve it is inscribed in.  Therefore $\kappa(K)\geq 4\pi$. To 
get a strict inequality, observe that a knot is not coplanar. By
repeatedly adding vertices to an alternating quadrisecant, the newly formed inscribed polygon eventually has four vertices non-coplanar, giving
it, and hence the knot, total curvature strictly greater than
$4\pi$. (See \lem{curvature} from \cite{Mil} below.)
\endproof

\begin{lemma}\label{lem:curvature}
Adding a new vertex to a closed polygon cannot decrease its total
curvature. The total curvature must increase if the new vertex $a_j$ and
three adjacent vertices ($a_{j-1}$, $a_{j+1}$, $a_{j+2}$) are not coplanar. \qed
\end{lemma}

The convex hull of a connected set $K$ in $\mathbb{R}^3$ is
characterized by the fact that every plane through a point in the hull
must intersect $K$. If $K$ is a closed curve, then a generic plane
must intersect $K$ an even number of times. Thus every plane through
each point of the convex hull is cut by $K$ at least twice.

This idea may be generalized as in \cite{ckks}:

\begin{definition}\label{def:hull} Let $K$ be a closed curve in
  $\mathbb{R}^3$. Its 
\emph{$n$th hull} $h_n(K)$ is the set of points $p\in\mathbb{R}^3$ such that
$K$ cuts every plane $P$ through $p$ at least $2n$-times.
If the intersections are transverse, then (thinking of $P$
as horizontal, and orienting $K$) there are equal numbers of upward
and downward intersections.
To handle non-transverse intersections, we again orient~$K$ and adopt
the following conventions. First, if $K\subset P$, we say $K$ cuts $P$
twice (once in either direction). If $K\cap P$ has infinitely many
components, then we say $K$ cuts $P$ infinitely often. Otherwise, each
connected component of the intersection is preceded and followed by
open arcs in $K$, with each lying to one side of $P$. An \emph{upward
intersection} will mean a component of $K\cap P$ preceded by an arc
below $P$ \emph{or} followed by an arc above $P$. (Similarly, a
\emph{downward intersection} will mean a component preceded by an arc
above $P$ or followed by an arc below $P$.) A glancing intersection,
preceded and followed by arcs on the same side of $P$, thus counts
twice, as both an upward and a downward intersection.
\end{definition}

Alternating quadrisecants give the geometric intuition that a knot
must travel twice around 
some point in space. Milnor observed that for a knotted curve, there
are planes in every direction which cut the knot four times. More
generally, there are points through which every plane cuts the knots
four times. But these points are precisely the second hull. In other words:

\begin{theorem}\label{thm:2hull} A knotted curve has
nonempty second hull.
\end{theorem}

This result was originally proved in \cite{ckks}. This paper
conjectured the existence of alternating quadrisecants as another way
of proving that the second hull of knotted curves is nonempty.   

\proof  By \thm{mainthma} we know that a knotted curve $K$ has an
alternating quadrisecant $abcd$. Let $t$ be a point on the midsegment $\overline{bc}$. We
claim that $t$ is in the second hull.

Project the knot radially to the unit sphere about $t$. Let the points
$a,b$ be at the north pole ($N$) and $c,d$ be at the south pole ($S$) of the
sphere. As $abcd$ is an alternating quadrisecant we see that the
projected knot visits the poles in the order $NSNS$. To show that $t$
is in the second hull, it suffices to show 
that the knot projection intersects any great circle at least four
times. The same conventions given in \defn{hull} apply to counting
intersections of $K$ with great circles.

There are two cases. Either the great circles are meridians (passing
through the north and south pole) or they are not. Suppose the great
circles are not meridians. Then as the projected curve visits the
poles in order $NSNS$, it must cut such a great circle at least four
times. Now suppose the great circle is a meridian. We divide the
projected knot into four arcs, each arc is the part of~$K$ between a 
north and south pole. If all of these arcs lie on the meridian then
this counts as an infinite number of intersections.  If one of these
arcs deviates from the meridian, this counts as two
intersections. Thus if two or more arcs leave the meridian, then there
are at least four intersections. Suppose that only one arc deviates
from the meridian. Such a curve stays in the plane 
determined by the meridian most of the time and goes to one side
once. This gives a knot of bridge number 1, which by \cite{Mil} is the unknot.
Thus for a knotted curve, the projected curve intersects all great
circles at least four times and the second hull is non-empty.
\endproof

The final application uses the second Main Theorem about essential
alternating quadrisecants to improve the known lower bounds of the
ropelength of knotted curves from 12 to 15.66.  

\begin{definition} Ropelength is the (scale invariant) quotient of
  length over thickness.
\end{definition}

The thickness $\tau(K)$ of a space curve is defined \cite{gm} to be twice the infimal
radius $r(x,y,z)$ of a circle passing through any three distinct
points of $K$. 
J.~Cantarella, R.~Kusner and J.M.~Sullivan~\cite{cks} showed that $\tau(K)=0$
unless~$K$ is~$C^{1,1}$. If $K$ is $C^1$ we can define normal tubes
abound~$K$, and then $\tau(K)$ is the supremal radius such that this tube
remains embedded. Cantarella, Kusner and Sullivan also proved that any (tame) knot or link has a
ropelength minimizer and gave lower bounds for the ropelength of
links \cite{cks}. These are sharp for certain simple cases, where each component
of the link is planar. However, these examples are the only known
ropelength minimizers, it seems difficult to describe explicitly the
shape of any tight knot. 

Numerical experiments \cite{Pie, Sul} suggest that the minimum
ropelength for a trefoil is slightly less than 16.4 and that there is
another tight trefoil with different symmetry and ropelength 18.7. The
best lower bound in \cite{cks} was 10.726, this was improved by Diao
\cite{Diao} who showed that any knot has ropelength more than~12.

In \cite{dds} we use the idea of essential quadrisecants to get better lower
bounds for ropelength. For each of the three types of quadrisecants we
use geometric arguments to obtain a lower bound for the ropelength of
the knot having a quadrisecant of that type. The worst of these three
cases is~13.936. However, from \cor{altquadess} we know that
any nontrivial~$C^{1,1}$ knot has an essential alternating
quadrisecant. This result together with the ropelength bound for a
knot with an essential alternating quadrisecant shows that
any nontrivial knot has ropelength at least~15.66.

\section*{Acknowledgements}
The author would like to acknowledge J.M.~Sullivan for suggesting this problem and for his endless support and advice. Thanks also to S.B.~Alexander, R.L.~Bishop and J.~Canteralla for their corrections and many helpful suggestions. Finally thanks to G.~Francis for the use of two of his pictures and for his help with the other figures. This work was partially supported by a Bourgain Fellowship whilst the author was a graduate student at the University of Illinois at Urbana-Champaign.


\end{document}